
\magnification=1200
\vsize = 8.9truein
\hsize = 6.5truein
\hoffset = -.15truecm
\voffset = -.5truecm
\parindent=0cm   \parskip=0pt
\pageno=1
\def\ind{\par\hskip 1cm\relax}
\pretolerance=500 \tolerance=1000  \brokenpenalty=5000

\catcode`\@=11

\font\eightrm=cmr8         \font\eighti=cmmi8
\font\eightsy=cmsy8        \font\eightbf=cmbx8
\font\eighttt=cmtt8        \font\eightit=cmti8
\font\eightsl=cmsl8        \font\sixrm=cmr6
\font\sixi=cmmi6           \font\sixsy=cmsy6
\font\sixbf=cmbx6


\font\tengoth=eufm10       \font\tenbboard=msym10
\font\eightgoth=eufm8      \font\eightbboard=msym8
\font\sevengoth=eufm7      \font\sevenbboard=msym7
\font\sixgoth=eufm6        \font\fivegoth=eufm5

\skewchar\eighti='177 \skewchar\sixi='177
\skewchar\eightsy='60 \skewchar\sixsy='60


\newfam\gothfam           \newfam\bboardfam
\newfam\cyrfam

\def\tenpoint{%
  \textfont0=\tenrm \scriptfont0=\sevenrm \scriptscriptfont0=\fiverm
  \def\rm{\fam\z@\tenrm}%
  \textfont1=\teni  \scriptfont1=\seveni  \scriptscriptfont1=\fivei
  \def\oldstyle{\fam\@ne\teni}\let\old=\oldstyle
  \textfont2=\tensy \scriptfont2=\sevensy \scriptscriptfont2=\fivesy
  \textfont\gothfam=\tengoth \scriptfont\gothfam=\sevengoth
  \scriptscriptfont\gothfam=\fivegoth
  \def\goth{\fam\gothfam\tengoth}%
  \textfont\bboardfam=\tenbboard \scriptfont\bboardfam=\sevenbboard
  \scriptscriptfont\bboardfam=\sevenbboard
  \def\bb{\fam\bboardfam\tenbboard}%
  \textfont\itfam=\tenit
  \def\it{\fam\itfam\tenit}%
  \textfont\slfam=\tensl
  \def\sl{\fam\slfam\tensl}%
  \textfont\bffam=\tenbf \scriptfont\bffam=\sevenbf
  \scriptscriptfont\bffam=\fivebf
  \def\bf{\fam\bffam\tenbf}%
  \textfont\ttfam=\tentt
  \def\tt{\fam\ttfam\tentt}%
  \abovedisplayskip=9pt plus 3pt minus 9pt
  \belowdisplayskip=\abovedisplayskip
  \abovedisplayshortskip=0pt plus 3pt
  \belowdisplayshortskip=4pt plus 3pt
  \smallskipamount=3pt plus 1pt minus 1pt
  \medskipamount=6pt plus 2pt minus 2pt
  \bigskipamount=12pt plus 4pt minus 4pt
  \normalbaselineskip=12pt
  \setbox\strutbox=\hbox{\vrule height8.5pt depth3.5pt width0pt}%
  \let\bigf@nt=\tenrm       \let\smallf@nt=\sevenrm
  \normalbaselines\rm}

\def\eightpoint{%
  \textfont0=\eightrm \scriptfont0=\sixrm \scriptscriptfont0=\fiverm
  \def\rm{\fam\z@\eightrm}%
  \textfont1=\eighti  \scriptfont1=\sixi  \scriptscriptfont1=\fivei
  \def\oldstyle{\fam\@ne\eighti}\let\old=\oldstyle
  \textfont2=\eightsy \scriptfont2=\sixsy \scriptscriptfont2=\fivesy
  \textfont\gothfam=\eightgoth \scriptfont\gothfam=\sixgoth
  \scriptscriptfont\gothfam=\fivegoth
  \def\goth{\fam\gothfam\eightgoth}%
    \textfont\bboardfam=\eightbboard \scriptfont\bboardfam=\sevenbboard
  \scriptscriptfont\bboardfam=\sevenbboard
  \def\bb{\fam\bboardfam}%
  \textfont\itfam=\eightit
  \def\it{\fam\itfam\eightit}%
  \textfont\slfam=\eightsl
  \def\sl{\fam\slfam\eightsl}%
  \textfont\bffam=\eightbf \scriptfont\bffam=\sixbf
  \scriptscriptfont\bffam=\fivebf
  \def\bf{\fam\bffam\eightbf}%
  \textfont\ttfam=\eighttt
  \def\tt{\fam\ttfam\eighttt}%
  \abovedisplayskip=9pt plus 3pt minus 9pt
  \belowdisplayskip=\abovedisplayskip
  \abovedisplayshortskip=0pt plus 3pt
  \belowdisplayshortskip=3pt plus 3pt
  \smallskipamount=2pt plus 1pt minus 1pt
  \medskipamount=4pt plus 2pt minus 1pt
  \bigskipamount=9pt plus 3pt minus 3pt
  \normalbaselineskip=9pt
  \setbox\strutbox=\hbox{\vrule height7pt depth2pt width0pt}%
  \let\bigf@nt=\eightrm     \let\smallf@nt=\sixrm
  \normalbaselines\rm}
\tenpoint
\def\pc#1{\bigf@nt#1\smallf@nt}         \def\pd#1 {{\pc#1} }
\def\up#1{\raise 1ex\hbox{\smallf@nt#1}}
\catcode`\@=12
\mathcode`A="7041 \mathcode`B="7042 \mathcode`C="7043 \mathcode`D="7044
\mathcode`E="7045 \mathcode`F="7046 \mathcode`G="7047 \mathcode`H="7048
\mathcode`I="7049 \mathcode`J="704A \mathcode`K="704B \mathcode`L="704C
\mathcode`M="704D \mathcode`N="704E \mathcode`O="704F \mathcode`P="7050
\mathcode`Q="7051 \mathcode`R="7052 \mathcode`S="7053 \mathcode`T="7054
\mathcode`U="7055 \mathcode`V="7056 \mathcode`W="7057 \mathcode`X="7058
\mathcode`Y="7059 \mathcode`Z="705A

\def\spacedmath#1{\def\packedmath##1${\bgroup\mathsurround=0pt ##1\egroup$}%
\mathsurround#1 \everymath={\packedmath}\everydisplay={\mathsurround=0pt }}
\def\nospacedmath{\mathsurround=0pt \everymath={}\everydisplay={} }
\def\diagram#1{\def\normalbaselines{\baselineskip=0truept
\lineskip=10truept\lineskiplimit=1truept}   \matrix{#1}}

\def\hfl#1#2{\normalbaselines{\baselineskip=0truept
\lineskip=10truept\lineskiplimit=1truept}\nospacedmath\smash\ {\mathop
{\hbox to 12truemm{\rightarrowfill}\
}\limits^{\scriptstyle#1}_{\scriptstyle#2}}}

\def\ghfl#1#2{\normalbaselines{\baselineskip=0pt
\lineskip=10truept\lineskiplimit=1truept}\nospacedmath\smash{\mathop{\hbox to
16truemm{\rightarrowfill}}\limits^{\scriptstyle#1}_{\scriptstyle#2}}}

\def\phfl#1#2{\normalbaselines{\baselineskip=0pt
\lineskip=10truept\lineskiplimit=1truept}\nospacedmath\smash\
{\mathop{\hbox to 8truemm{\rightarrowfill}\
}\limits^{\scriptstyle#1}_{\scriptstyle#2}}}

\def\vfl#1#2{\llap{$\scriptstyle#1$}\left\downarrow\vbox to
6truemm{}\right.\rlap{$\scriptstyle#2$}}
 \def\cqfd{\kern 2truemm\unskip\penalty 500\vrule
height 4pt depth 0pt width 4pt\medbreak} \def\carre{\vrule height 4pt depth
0pt width 4pt} \def\moins{\mathrel{\hbox{\vrule height 3pt depth -2pt
width 6pt}}}
\def\rond{\kern 1pt{\scriptstyle\circ}\kern 1pt}
\def\epi{\rightarrow \kern -3mm\rightarrow }
\def\longepi{\longrightarrow \kern -5mm\longrightarrow }
\def\mono{\lhook\joinrel\mathrel{\longrightarrow}}
\def\End{\mathop{\rm End}\nolimits}
\def\Hom{\mathop{\rm Hom}\nolimits}
\def\Aut{\mathop{\rm Aut}\nolimits}
\def\im{\mathop{\rm Im}\nolimits}
\def\Ker{\mathop{\rm Ker}\nolimits}
\def\Coker{\mathop{\rm Coker}}
\def\det{\mathop{\rm det}\nolimits}
\def\Pic{\mathop{\rm Pic}\nolimits}

\def\dim{\mathop{\rm dim}\nolimits}

\def\Tr{\mathop{\rm Tr}\nolimits}
\def\bk{\backslash}
\def\Sp{\mathop{\rm Spec\,}\nolimits}
\def\iso{\mathrel{\mathop{\kern 0pt\longrightarrow }\limits^{\sim}}}
\def\vb{vector bundle }
\def\al{$k$\kern -1.5pt -algebr}
\def\rk{\mathop{\rm rk\,}\nolimits}
\def\div{\mathop{\rm div\,}\nolimits}
\def\Ad{\mathop{\rm Ad}\nolimits}
\def\Res{\mathop{\rm Res}\nolimits}
\def\Lie{\mathop{\rm Lie}\nolimits}
\def\limind{\mathop{\oalign{lim\cr\hidewidth$\longrightarrow
$\hidewidth\cr}}} \def\Ann{\mathop{\rm Ann}\nolimits}

\def\codim{\mathop{\rm codim}\nolimits}
\baselineskip 15pt
\spacedmath{2pt}
\overfullrule=0pt

\font\gros=cmr10 scaled\magstep3 
\null
\centerline{{\bf \gros Conformal blocks and generalized theta functions}}
\medskip

\centerline{Arnaud {\pc BEAUVILLE} and Yves {\pc LASZLO}
\footnote{\parindent 0.5cm(*)}{\eightpoint Both authors were partially
supported by the European Science Project ``Geometry of Algebraic Varieties",
Contract no. SCI-0398-C(A).}} \vskip1.5cm
\parindent=0cm
{\bf Introduction}
\ind The aim of this paper is to construct a canonical isomorphism between
two vector spaces associated to a Riemann surface $X$. The first of these
spaces is the space of {\it conformal blocks}  $B_c(r)$ (also called the
space of vacua), which plays an important role in conformal field theory. It
is defined as follows: choose a point $p\in X$, and let $A_X$ be the ring of
algebraic functions on $X\moins p$. To each integer $c\ge 0$ is associated a
representation $V_c$ of the Lie algebra ${\goth sl}_r\bigl({\bf
C}((z))\bigr)$, the {\it basic representation} of level $c$ (more correctly
it is a representation of the universal extension of ${\goth sl}_r\bigl({\bf
C}((z))\bigr)$ -- see \S 7 for details). The ring $A_X$ embeds into ${\bf
C}((z))$ by associating to a function its Laurent development at $p$; then
$B_c(r)$ is the space of linear forms on $V_c$ which vanish on the elements
$A(z)v$ for $A(z)\in {\goth sl}_r(A_X)\ ,\ v\in V_c$. \ind The second space
comes from algebraic geometry, and is defined as follows. Let ${\cal
SU}_X(r)$ be the moduli space of semi-stable  rank $r$ vector bundles on $X$
with trivial determinant.  One can define a theta divisor on ${\cal SU}_X(r)$
in the same way one does in the rank $1$ case: one chooses a line bundle $L$
on $X$ of degree $g-1$, and considers the locus of vector bundles $E\in{\cal
SU}_X(r)$ such that $E\otimes L$ has a nonzero section. The associated line
bundle ${\cal L}$ is called the {\it determinant bundle}; the space we are
interested in  is $H^0({\cal SU}_X(r),{\cal L}^c)$.  This space  can be
considered as a non-Abelian version of the space of $c^{\rm th}$-order theta
functions on the Jacobian of $X$, and is  sometimes called the space of
generalized theta functions. We will prove that it is canonically isomorphic
to $B_c(r)$. By [T-U-Y] this implies that the space $H^0({\cal SU}_X(r),{\cal
L}^c)$ satisfies the so-called {\it fusion rules}, which allow to compute its
dimension in a purely combinatorial way. According to a conjecture of
Verlinde, proved in [M-S], this dimension is given  by the famous Verlinde
formula ([V], see cor. 8.6).
 \ind The  isomorphism $B_c(r)\iso H^0({\cal SU}_X(r),{\cal L}^c)$ is
certainly known to the physicists -- see e.g. [W]. Our point is that this
can be proved in a purely mathematical way. In fact we hope to convince the
reader that even in an infinite-dimensional context, the methods of algebraic
geometry provide a flexible and efficient language (though a little
frightening at first glance!). \ind Our strategy is as follows. First, by
trivializing vector bundles on $X\moins p$ and on a neighborhood of $p$, we
construct a bijective correspondence between the moduli space and the double
coset space $SL_r(A_X)\bk SL_r\bigl({\bf C}((z))\bigr)/SL_r({\bf C}[[z]])$
(this is a quite classical idea which  goes back to A. Weil). Sections 1 to 3
are devoted to make sense of this as an isomorphism between geometric
objects. We show that the quotient ${\cal Q}:=SL_r\bigl({\bf C}((z))\bigr)/
SL_r({\bf C}[[z]])$ as well as  the group $SL_r(A_X)$ is an {\it
ind-variety}, that is a direct limit of an increasing sequence of algebraic
varieties. The quotient $SL_r(A_X)\bk {\cal Q}$ makes sense as a {\it stack}
(not far from what topologists call an orbifold), and  this stack is
canonically isomorphic to the moduli stack ${\cal SL}_X(r)$ of vector bundles
on $X$ with trivial determinant.

\ind The determinant line bundle ${\cal L}$ lives naturally on the moduli
stack, and the next step is to identify its pull back to ${\cal Q}$. In order
to do this we first construct the central ${\bf C}^*$-extension
$\widehat{SL}_r\bigl({\bf C}((z))\bigr)$ and the $\tau$ function on this
group, and show that the $\tau$ function defines a section of a line bundle
${\cal L}_\chi$ on ${\cal Q}$ (\S 4). We then prove that the pull back of
${\cal L}$ to ${\cal Q}$ is isomorphic to ${\cal L}_\chi$ (\S 5). A theorem
of Kumar and Mathieu  identifies the space $H^0({\cal Q},{\cal L}_{\chi}^c)$
with the dual $V_c^*$ of the basic representation $V_c$; it follows, almost
by definition of a stack, that $H^0({\cal SL}_X(r),{\cal L}^c)$ can be
identified with the elements of $V_c^*$ which are invariant under the group
$SL_r(A_X)$ (\S 7). This turns out to be the same as the linear forms
annihilated by the Lie algebra: the key point is that the group  $SL_r(A_X)$
is {\it reduced} (\S 6) -- a highly non-trivial property in our
infinite-dimensional set-up.  The final step is to prove that the sections of
${\cal L}^c$ on the moduli stack and on the moduli space are the same (\S 8)
-- this is essentially Hartogs' theorem, since the substack of non-stable
bundles is of codimension $\ge 2$. \ind In the last section we state and
prove the corresponding result for the moduli space of vector bundles of rank
$r$ and  determinant $L$ for any line bundle $L$ on $X$. \smallskip \ind The
methods of this paper should extend to the general case of principal bundles
under a semi-simple algebraic group $G$. We have chosen to work in the
context of vector bundles ({\it i.e. $G=SL_r({\bf C})$}) because it is  by
far the most important case for algebraic geometers, and it is easier to
explain in so far as it appeals very little to the rather technical machinery
of Kac-Moody groups. Also the general case can be to a large extent reduced
to this one. \smallskip

\ind Most of this work was done in the Spring of
1992, and we have  lectured in various places about it. In July 1992 we
heard of G.~Faltings beautiful ideas, which should prove at the same time
both our result and that of [T-U-Y] (in the more general case of principal
bundles). These ideas are sketched in [F], but (certainly due to our own
incompetence) we were unable to understand some of the key points in the
proof. We have therefore decided after some time to write a complete version
of  our proof, if only to provide an
 introduction to Faltings' ideas.

\ind Part of our results have been  obtained  independently (also in the
 context of principal bundles) by S. Kumar, Narasimhan and Ramanathan
[K-N-R].
\bigskip
{\sevenrm \ind We would like to thank G. Faltings for useful letters, and V.
Drinfeld for pointing out an inaccuracy in an earlier version of this paper.}
\vskip 1.7cm
{\bf 1. The ind-groups ${\bf GL}_r(K)$ and ${\bf SL}_r(K)$.}
\smallskip
{\it $k$-spaces and ind-schemes}
\ind (1.1) Throughout this paper we'll work over an algebraically closed
field  $k$ of characteristic $0$. A \al a will always be assumed to be
associative, commutative and unitary. Our basic objects will be {\it
$k$-spaces} in the sense of [L-MB]: by definition, a $k$-space (resp. a
$k$-group) is a functor $F$ from the category of \al as to the category of
sets (resp. of groups) which is a sheaf for the faithfully flat topology.
 Recall that this means that for any faithfully flat homomorphism
 $R\rightarrow R'$,  
the diagram $$F(R)\longrightarrow F(R') \raise -3pt\hbox{$\
\buildrel\displaystyle \longrightarrow\over{\longrightarrow}\ $}F(R'\otimes_R
R')$$  is exact; in most cases the verification that this is indeed the case
is quite easy, and will be left to the reader. The category of  schemes over
$k$ can be (and will be) viewed as a full subcategory of the category of
$k$-spaces. A scheme will always be assumed to be quasi-compact and
quasi-separated.

\ind Direct limits exist in the category of $k$-spaces; we'll say that a
$k$-space (resp. a $k$-group) is an {\it ind-scheme}
(resp. an {\it ind-group})
 if it is the direct limit of a directed system of schemes.
 Let $(X_\alpha)_{\alpha\in I}$ be a directed system of schemes,  $X$ its
limit in the category of $k$-spaces, and $S$ a $k$-scheme. The set ${\rm
Mor}\,(S,X)$ of morphisms of $S$ into $X$ is the direct limit of the sets
${\rm Mor}\,(S,X_\alpha)$, while the set ${\rm Mor}\,(X,S)$ is the inverse
limit of the sets ${\rm Mor}\,(X_\alpha,S)$. \bigskip

{\it The groups ${\bf GL}_r(K)$ and ${\bf GL}_r({\cal O})$}
\ind (1.2) Let $z$ be an indeterminate. We will denote by ${\cal O}$ the
formal series ring $k[[z]]$ and by $K$ the field $k((z))$ of meromorphic
formal series in $z$. We let ${\bf GL}_r({\cal O})$ (or ${\bf GL}_r(k[[z]])$)
be the $k$-group  $R\mapsto
 GL_r(R[[z]])$, and ${\bf GL}_r(K)$ (or ${\bf GL}_r\bigl(k((z))\bigr)$) be
the $k$-group  $R\mapsto GL_r\bigl(R((z))\bigr)$. We define in the same way
the $k$-groups ${\bf SL}_r({\cal O})$ and ${\bf SL}_r(K)$. For $N\ge 0$, we
denote by $G^{(N)}(R)$ (resp. $S^{(N)}(R)$) the set of matrices $A(z)$ in
$GL_r\bigl(R((z))\bigr)$ (resp. in $SL_r\bigl(R((z))\bigr)$) such that both
$A(z)$ and $A(z)^{-1}$ have  a pole of order $\le N$. This defines
subfunctors $G^{(N)}$ and $S^{(N)}$ of ${\bf GL}_r(K)$ and ${\bf SL}_r(K)$
respectively. \medskip
{\pc PROPOSITION} 1.2.-- {\it The $k$-group ${\bf GL}_r({\cal O})$} (resp.
${\bf SL}_r({\cal O})$) {\it is an affine group scheme. The $k$-group ${\bf
GL}_r(K)$} (resp. ${\bf SL}_r(K)$) {\it is an ind-group, direct limit of the
sequence of schemes $(G^{(N)})_{N\ge 0}$}  (resp.  $(S^{(N)})_{N\ge 0}$).

 \ind For any \al a $R$, let us denote by $M_r(R)$ the vector space of
$r$-by-$r$ matrices with entries in $R$. The set $GL_r(R[[z]])$ consists of
matrices $\displaystyle A(z)=\sum_{n\ge 0}A_n z^n$, with $A_0\in GL_r(R)$ and
$A_n\in M_r(R)$ for $n\ge 1$; therefore the group ${\bf GL}_r({\cal O})$ is
represented by the affine  scheme $\displaystyle GL_r(k)\times \prod_1^\infty
M_r(k)$.

\ind Let $M^{(N)}(R)$ be the
space of $r$-by-$r$ matrices $\displaystyle A(z)=\sum_{n\ge -N}A_n z^n$,
with $A_n\in M_r(R)$.  The functor $M^{(N)}$ is represented by the affine
scheme $ \displaystyle \prod_{n\ge -N}M_r(k) $, and the functor $G^{(N)}$ is
represented by a closed (affine)  subscheme of $M^{(N)}\times M^{(N)}$
(identify $G^{(N)}(R)$  with the subset of $M^{(N)}(R)\times M^{(N)}(R)$
consisting of couples $(A(z),B(z))$ such that $A(z)B(z) = I$).  One has
$\displaystyle GL_r\bigl(R((z))\bigr) = \bigcup_{N\ge 0}G^{(N)}(R)$, hence
{\it the $k$-group ${\bf GL}_r(K)$ is the
 direct limit of the sequence of schemes $(G^{(N)})_{N\ge 0}$.}

\ind Let $N$ be a non-negative integer. There exist universal polynomials
 \break $P_m^{(N)}\bigl((A_n)_{n\ge -N}\bigr)\quad (m\ge -rN)$ on the affine
space $\displaystyle \prod_{n\ge -N}M_r(k)$ such that the determinant of an
element  $\displaystyle A(z)=\sum_{n\ge -N}A_nz^n$
of $GL_r\bigl(R((z))\bigr)$ is given by  $$\det A(z) = \sum_{m\ge
-rN}P_m^{(N)}\bigl((A_n)_{n\ge -N}\bigr) z^m\quad .$$  It follows
that the functor $S^{(N)}$ is representable by a closed affine subscheme of
$G^{(N)}$. In particular, $S^{(0)}={\bf SL}_r({\cal O})$ is an affine scheme,
and ${\bf SL}_r(K)$ is an ind-scheme, direct limit of the sequence
$\bigl(S^{(N)}\bigr)_{N\ge 0}$.\cqfd
 \bigskip
{\it ${\bf GL}_r(K)$ and vector bundles}

\ind (1.3) We now start the geometric side of this paper; we fix once and
for all a smooth (connected) projective curve $X$ over $k$, and a closed
point $p$ of $X$. We put $X^*=X\moins p$. We denote by ${\cal O}$ the
completion of the local ring of $X$ at $p$, and by $K$ its field of
fractions. We will choose\footnote{\parindent .5cm $^1$}{\eightpoint see
Remark (1.7) below.} \parindent=0cm a local coordinate $z$ at $p$ and
identify ${\cal O}$ with $k[[z]]$ and $K$ with $k((z))$.
  Let $R$ be a  $k$-algebra. We put $X_R = X\times_k \Sp(R)$, $X_R^* =
 X^*\times_k \Sp(R)$, $D_R = \Sp\bigl(R[[z]]\bigr)$ and $D_R^* =
\Sp\bigl(R((z))\bigr)$. We consider the cartesian diagram \vskip -10pt
\def\hfl#1#2{\normalbaselines{\baselineskip=0truept
\lineskip=10truept\lineskiplimit=1truept}\nospacedmath\smash{\mathop{\hbox to
12truemm{\rightarrowfill}}\limits^{\scriptstyle#1}_{\scriptstyle#2}}}

$$\diagram{
D_R^*& \lhook\joinrel\mathrel{\hfl{}{}}&D_R&\cr
\vfl{}{}&&\vfl{}{f}&\cr
X_R^*&\lhook\joinrel\mathrel{\hfl{j}{}}&X_R&\ .\cr
}\leqno (1.3)$$\vskip -10pt
\spacedmath{2pt}
\ind When $R=k$, we may think of $f(D)$ as a small disk in $X$ around $p$,
and of $f(D^*)$ as the punctured disk $f(D)\moins p$. We want to say that the
ind-group ${\bf GL}_r(K)$ parametrizes bundles which are trivialized on $X^*$
and on $D$.

\ind We consider triples $(E,\rho,\sigma)$ where $E$ is a vector bundle on
$X_R$,\break $\rho : {\cal O}_{X_R^*}^r \longrightarrow E_{|X_R^*}$ a
trivialization of $E$ over $X_R^*$, $\sigma : {\cal
O}^r_{D_R}\longrightarrow E_{|D_R}$ a
trivialization  of $E$ over $D_R$.  We let $T(R)$ be the set of isomorphism
 classes of triples $(E,\rho,\sigma)$ (with the obvious notion of
isomorphism). \bigskip

{\pc PROPOSITION} 1.4. -- {\it The ind-group ${\bf GL}_r(K)$ represents the
functor $T$.} \smallskip

\ind Let $(E,\rho,\sigma)$ be an element of $T(R)$. Pulling back the
trivializations $\rho$ and $\sigma$ to $D^*_R$ provides two trivializations
$\rho^*$ and $\sigma^*$ of the pull back of $E$ over ${D^*_R}$: these
trivializations differ by an element $\gamma = \rho^{*-1}\rond\sigma^*$ of $
GL_r\bigl(R((z))\bigr)$.

\def\wfl#1#2{\llap{$\scriptstyle#1$}\left\uparrow\vbox to 6truemm{}\right.
\rlap{$\scriptstyle#2$}}

\ind Let us now drop the suffix
$R$ to simplify the notation. Let  ${\cal K}_D$ be the quasi-coherent sheaf
on $D$ associated to the $R[[z]]$-module $R((z))$.  We have an exact sequence
$$0\rightarrow {\cal O}_X\longrightarrow j_*{\cal
O}_{X^*}\longrightarrow f_*({\cal K}_D/{\cal O}_D)\rightarrow 0\ .$$
Tensoring with $E$ and using the trivializations $\rho$ and $\sigma$, we get a
commutative diagram with exact rows $$\diagram{0\rightarrow
&E&\longrightarrow &j_*(E_{|X^*})&\longrightarrow
&f_*\Bigl(E_{|D}\otimes{\cal K}_D/{\cal O}_D)&\rightarrow 0\cr
&\left\Vert\vbox to 6truemm{}\right.&&\wfl{}{\rho}&&\wfl{}{\sigma}\cr
0\rightarrow &E&\longrightarrow &j_*{\cal O}_{X^*}^r&\phfl{\overline\gamma}{}
&f_*({\cal K}_D/{\cal O}_D)^r&\rightarrow 0\cr}$$ where $\overline{\gamma}$
is the composition of the natural map $j_*{\cal O}_{X^*}^r\longrightarrow
f_*({\cal K}_D)$, the automorphism $\gamma^{-1}$ of $f_*({\cal K}_D)^r$, and
the canonical projection $f_*({\cal K}_D)^r\longrightarrow f_*({\cal
K}_D/{\cal O}_D)^r$. \ind Conversely, let us start from an element $\gamma$
of $GL_r\bigl(R((z))\bigr)$. We claim that the homomorphism
$\overline{\gamma}:j_*{\cal O}_{X^*}^r\longrightarrow  f_*({\cal K}_D/{\cal
O}_D)^r$ defined by the above recipe is surjective, and that its kernel
$E_\gamma$ is locally free of rank $r$. By descent theory it is enough to
check these assertions after pull back to $X^*$ and to $D$. They are clear
over $X^*$, since the exact sequence reduces to an isomorphism
$\rho^{-1}:E_\gamma\rightarrow {\cal O}_{X^*}^r$. Over $D$, we observe that
the canonical map $f^*f_*({\cal K}_D/{\cal O}_D)\longrightarrow {\cal
K}_D/{\cal O}_D$ is an isomorphism (express for instance ${\cal K}_D/{\cal
O}_D$ as the  limit of the direct system $(\cdots\rightarrow {\cal
O}_D/(z^n)\phfl{z}{}{\cal O}_D/(z^{n+1})\longrightarrow \cdots)$). Therefore
we get an exact sequence $$0\rightarrow f^*E\longrightarrow {\cal
K}_D^r\hfl{p\rond\gamma^{-1}}{}({\cal K}_D/{\cal O}_D))^r\rightarrow 0\
,$$where $p:{\cal K}_D^r\longrightarrow ({\cal K}_D/{\cal O}_D))^r$ is the
canonical map. In other words, $\gamma$ induces an isomorphism $\sigma:{\cal
O}_D^r\rightarrow f^*E$. Thus $E_\gamma$ is a vector bundle, so we have
associated to $\gamma$ a triple $(E_\gamma,\rho,\sigma)$ in $T(R)$. The two
constructions are clearly inverse of each other, hence the proposition. \cqfd
\bigskip
\ind From this proposition we get immediately
\medskip
 {\pc PROPOSITION} 1.5. -- {\it The ind-group ${\bf SL}_r(K)$ represents the
subfunctor $T_0$ of $T$ which associates to a \al a $R$ the set of
isomorphism classes of triples $(E,\rho,\sigma)$ where $E$ is a vector bundle
on $X_R$, $\rho : {\cal O}_{X_R^*}^r\longrightarrow E_{|X_R^*}$ and $\sigma
:  {\cal O}^r_{D_R}\longrightarrow E_{|D_R}$ are isomorphisms such that
$\wedge^r\rho$ and $\wedge^r\sigma$ coincide over $D^*_R$.}\cqfd
 \bigskip
{\it Remark} 1.6.-- The condition that the trivializations $\wedge^r\rho$ and
$\wedge^r\sigma$ coincide over $D^*_R$ means that they come from a global
trivialization of $\bigwedge^r E$. So we can rephrase prop. 1.5 by saying
that $T_0(R)$ is the  set of isomorphism classes of data
$(E,\rho,\sigma,\delta)$ where $\delta$ is a trivialization of $\bigwedge^r
E$, $\rho$ and $\sigma$ are trivializations of $E_{|X_R^*}$ and $E_{|D_R}$
respectively, such that $\wedge^r\rho$  coincide with $\delta_{|X_R^*}$
and $\wedge^r\sigma$ with $\delta_{|D_R}$. \smallskip
 {\it Remark} 1.7.-- The
R-algebras $R[[z]]$ and $R((z))$ do not actually depend on the choice of a
local coordinate $z$ at $p$:  $R[[z]]$ is the completion of the tensor
product $R\otimes_k{\cal O}$ with respect to the $(R\otimes {\goth m})$-adic
topology, where ${\goth m}$ is the maximal ideal of ${\cal O}$, and $R((z))$
is
 $R[[z]]\otimes_{\cal O}K$.

\bigskip

\ind Let us specialize prop. 1.4 and 1.5 to the case  $R=k$:
\smallskip

{\pc COROLLARY} 1.8. -- {\it Let us denote by $A_X$ the affine algebra
$\Gamma(X\moins p\, ,{\cal O}_X)$. There is a canonical bijective
correspondence between the set of isomorphism classes of rank $r$ vector
bundles on $X$ with trivial determinant {\rm  (resp.} with determinant of the
form ${\cal O}_X(n\,p)$ for some integer $n${\rm )} and the double coset
space $SL_r(A_X)\backslash SL_r(K)/SL_r({\cal O})$ } (resp.
$GL_r(A_X)\backslash GL_r(K)/GL_r({\cal O})$).
 \smallskip
 \ind Since two trivializations of $E_{|D}$ differ by an element of $
GL_r({\cal O})$, and two trivializations of  $E_{|X^*}$  by an element of $
GL_r(A_X)$, we deduce from prop. 1.4 a bijection between $GL_r(A_X)\backslash
GL_r(K)/GL_r({\cal O})$ and the set of isomorphism classes of rank $r$ vector
bundles  on $X$ which are trivial on $X^*$. But a projective module over a
Dedekind ring is free if and only if its determinant is free ([B], ch. 7, \S
4, prop.~24), hence our assertion for $GL_r$. The same proof applies for $
SL_r$.\cqfd  \bigskip
{\it Remark} 1.9.-- Let $\gamma$ be an element of
$GL_r(K)$, and $E$ the corresponding vector bundle.   Then the element
$\det\gamma$ of $GL_1(K)=K^*$ corresponds to the line bundle
$\det E\cong{\cal O}(n\,p)$. It
follows that {\it the integer $n$ such that $\det E\cong{\cal O}(n\,p)$ is the
order of
the Laurent series $\det\gamma$.}  This shows that the group ${\bf GL}_r(K)$
is the union of countably many open and closed sub-ind-schemes, distinguished
by the order of the determinant. \bigskip

\ind (1.10) Our first goal in the following sections will be to show that the
bijection defined in cor. 1.8 comes actually from an isomorphism between
algebro-geometric objects. Let us observe here that  the functor $R\mapsto
SL_r(A_{X_R})$ is a $k$-group, which will play an important r\^ole in our
story; we denote it by   ${\bf SL}_r(A_X)$. It is actually an ind-variety,
 limit of the affine varieties $\Gamma^{(N)}$ parametrizing matrices
$A=(a_{ij})$ with $\det A=1$
and $a_{ij}\in H^0(X,{\cal O}_X(Np))$ for all $i,j$.
We shall study this group in more detail in \S 7.
\bigskip

{\it Application: the Birkhoff decomposition}
\ind Let us apply  cor. 1.8 when $X={\bf P}^1$, and
$p=0$. The vector bundles on ${\bf P}^1$ with rank $r$ and trivial
determinant are parametrized by sequences of integers ${\bf
d}=(d_1,\ldots,d_r)$ with $d_1\le \ldots\le d_r$ and $\sum d_i=0$: to such a
sequence corresponds the vector bundle ${\cal O}_{{\bf
P}^1}(d_1)\oplus\ldots\oplus{\cal O}_{{\bf P}^1}(d_r)$, which is defined by
the diagonal matrix  $z^{\bf d}:={\rm diag}\,(z^{d_1},\ldots,z^{d_r})$. The
\al a $A_{{\bf P}^1}$ is simply $k[z^{-1}]$. We obtain the {\it Birkhoff
decomposition} $$SL_r(K)=\bigcup_{\bf d}\ SL_r(k[z^{-1}]) \ z^{\bf d}\
SL_r({\cal O})\ .\leqno (1.11)$$

\ind We shall see that the {\it big cell} ${\bf SL}_r(K)^0:={\bf
SL}_r(k[z^{-1}])\ {\bf SL}_r({\cal O})$ is open in ${\bf SL}_r(K)$. More
precisely, let us denote (abusively) by ${\bf SL}_r({\cal O}_-)$ the (closed)
sub-ind-group of ${\bf SL}_r(k[z^{-1}])$ parametrizing  matrices of the form
$A(z)=I+\sum_{n\ge 1}A_nz^{-n}$.  \medskip

 {\pc PROPOSITION} 1.12. -- {\it The multiplication map
 $\mu:{\bf SL}_r({\cal O}_-)\times {\bf SL}_r({\cal O})  \longrightarrow {\bf
SL}_r(K)$ is an open immersion.}
 \smallskip

\ind Let first $S$ be a scheme and  ${\cal E}$ a vector bundle over $S\times
{\bf P}^1$; we denote by    $q:S\times {\bf P}^1\rightarrow S$  the
projection map.  Let $S^0$ be the biggest open subset of $S$ over which the
canonical map $q^*q_*{\cal E}\rightarrow {\cal E}$ is an isomorphism; this is
the locus of points $s$ in $S$ such that ${\cal E}_{|\{s\}\times {\bf P}^1}$
is trivial\footnote{\parindent 0.4cm$^1$}{\vtop{\eightpoint\baselineskip 12pt
\noindent We are using here (and will use in the sequel) the fact that the
standard base change theorems for  coherent cohomology are valid without any
noetherian hypothesis for projective morphisms (see [SGA 6], exp. III).}}. If
moreover the bundle ${\cal E}_{|S\times\{0\}}$ is trivial, so is the
restriction of ${\cal E}$ to $S^0$.

\ind We apply these remarks to $S={\bf SL}_r(K)$. Let $R$ be any
\al a. Clearly the map ${\bf SL}_r(z^{-1}R[z^{-1}])\times  {\bf
SL}_r(R[[z]]) \longrightarrow {\bf SL}_r(R((z)))$ is injective, and its image
corresponds to triples $(E,\rho,\sigma)$ over ${\bf P}^1_R$ where the vector
bundle $E$ is  trivial. Therefore $\mu$ induces an isomorphism from ${\bf
SL}_r(z^{-1}k[z^{-1}])\times{\bf SL}_r({\cal O}) $ onto the open
sub-ind-scheme ${\bf SL}_r(K)^0$.\cqfd

\vskip 1.7cm
{\bf 2. The homogeneous space $ {\bf SL}_r(K)/{\bf SL}_r({\cal O})$}
\smallskip
\ind In the preceding section we have described (cor. 1.8) a bijection
between the set of isomorphism classes of rank $r$ vector bundles on $X$ with
trivial determinant and the double coset space $SL_r(A_X)\backslash
SL_r(K)/SL_r({\cal O})$. Our aim in this section and the following is to show
that this gives in fact a description  of the moduli space -- actually of the
moduli stack.
We therefore  need to understand the algebraic structure of the set
$SL_r(A_X)\backslash SL_r(K)/SL_r({\cal O})$. We'll start with the quotient
${\bf SL}_r(K)/{\bf SL}_r({\cal O})$, which will turn out to be as nice as we
 can reasonably hope, namely a direct limit of projective varieties (thm. 2.5
below).    Let us first recall that such a quotient always exists as a
$k$-space -- it is simply the sheaf (for the faithfully flat topology)
associated to the presheaf $R\mapsto  SL_r\bigl(R((z))\bigr)/SL_r(R[[z]])$.
 \medskip

\bigskip
{\pc PROPOSITION} 2.1.-- {\it The $k$-space ${\cal Q}:={\bf SL}_r(K)/{\bf
SL}_r({\cal O})$ represents the functor which associates to a $k$-algebra $R$
the set of isomorphism classes of pairs $(E,\rho)$, where $E$ is a vector
bundle over $X_R$  and $\rho$ a trivialization of $E$  over $X^*_R$ such that
$\wedge^r\rho$ extends to a trivialization of $\bigwedge^r E$.}
 \smallskip
\ind This is just standard descent theory: let $R$ be a \al a and $q$  an
element of ${\cal Q}(R)$. By definition there exists a faithfully flat
homomorphism $R\rightarrow R'$ and an element $\gamma$ of
$SL_r\bigl(R'((z))\bigr)$ such that the image of $q$ in ${\cal Q}(R')$ is the
class of $\gamma$.  To $\gamma$ corresponds by prop. 1.5  a triple
$(E',\rho',\sigma')$ over $X_{R'}$. Let $R''=R'\otimes_R R'$, and let
$(E''_1,\rho''_1)$, $(E''_2,\rho''_2)$ denote the pull-backs of $(E',\rho')$
by the two projections of $X_{R''}$ onto $X_{R'}$. Since the two images of
$\gamma$ in $SL_r\bigl(R''((z))\bigr)$ differ by an element of
$SL_r\bigl(R''[[z]]\bigr)$, these pairs  are isomorphic; this means that the
isomorphism $\rho''_2 \rho_1''^{-1}$ over $X_{R''}^*$ extends to an
isomorphism $u:E''_1\rightarrow E''_2$ over $X_{R''}$. This isomorphism
satisfies the usual cocycle condition, because it is enough to check it over
$X^*$, where it is obvious. Therefore $(E',\rho')$ descends to a pair
$(E,\rho)$ on $X_R$ as in the statement of the proposition.

\ind Conversely, given a pair $(E,\rho)$ as above over $X_R$, we can find a
faithfully flat homomorphism $R\rightarrow R'$ and a trivialization $\sigma'$
of the pull back
of $E$ over $D_{R'}$ such that $\wedge^r\sigma'$ coincides with
$\wedge^r\rho$ over $D^*_{R'}$ (in fact  $\Sp(R)$ is covered by open subsets
$\Sp(R_\alpha)$ such that $E$ is trivial over $D_{R_\alpha}$, and we can take
$R'=\prod R_\alpha$). By prop. 1.5 we get an element $\gamma'$ of
$SL_r\bigl(R'((z))\bigr)$ such that the two images of $\gamma'$ in
$SL_r\bigl(R''((z))\bigr)$ (with $R''=R'\otimes_R R'$) differ by an element
of $SL_r\bigl(R''[[z]]\bigr)$; this gives an element of ${\cal Q}(R)$. The
two constructions are clearly inverse one of each other.\cqfd
\bigskip
{\it Remark} 2.2.-- The same construction applies without any modification
to the case of the group $GL_r$, giving a $k$-space ${\bf GL}_r(K)/{\bf
GL}_r({\cal O})$ which represents the functor associating to a  \al a $R$ the
set of isomorphism classes of pairs $(E,\rho)$, where $E$ is a vector bundle
over $X_R$  and $\rho$ a trivialization of $E$.

\ind Given such a pair, $\wedge^r \rho$ is a rational section of
$\bigwedge^r E$, and therefore extends to a trivialization of $(\bigwedge^r
E)(np)$ for some locally constant function $n:\Sp(R)\rightarrow {\bf Z}$.
Hence {\it the quotient ${\bf GL}_r(K)/{\bf GL}_r({\cal O})$ is a disjoint
union of a sequence $({\cal Q}_d)_{d\in{\bf Z}}$ of $k$-spaces}: the space
${\cal Q}_d$ parametrizes the pairs $(E,\rho)$ where $\wedge^r\rho$ extends
to an isomorphism ${\cal O}_{X_R}(dp)\iso \bigwedge^rE$. In group-theoretic
terms, ${\cal Q}_d$ is the quotient  ${\bf GL}_r(K)^{(d)}/{\bf GL}_r({\cal
O})$, where ${\bf GL}_r(K)^{(d)}$ is the open and closed sub-ind-scheme of
${\bf GL}_r(K)$ parametrizing matrices $A(z)$ such that the Laurent
series\footnote{\parindent 0.5cm $^1$}{\eightpoint A Laurent series
$\varphi\in R((z))$ is said of order $d$ if its image in $F((z))$ has order
$d$ for each homomorphism of $R$ into a field $F$.} \parindent=0cm
 $\det A(z)$ has order $d$ (1.9). If $\gamma_d$ is an element of
$GL_r(K)^{(d)}$,   left multiplication by $\gamma_d$
defines an isomorphism of ${\cal Q}$ onto ${\cal Q}_d$.

 \bigskip
{\it ${\cal Q}$ as a Grassmannian}

\ind Our quotient space ${\cal Q}$ is related to the infinite Grassmannian
used by the Japanese school (see [S-W]) in the following way. For any
$k$-algebra $R$, define a {\it lattice} in $R((z))^r$ as a
sub-$R[[z]]$-module $W$ of $R((z))^r$  which is projective of rank
$r$, and such that $\bigcup z^{-n}W=R((z))^r$.  It is an easy exercise in
Algebra to show that this amounts to say that $W$ is a sub-$R[[z]]$-module
of  $R((z))^r$  such that  $$z^N\,R[[z]]^r\i W \i z^{-N}\,R[[z]]^r$$
for some integer $N$, and such that the R-module $W/z^N\,R[[z]]^r$ is
projective. Let us say moreover that the lattice $W$ is {\it special} if the
lattice $\bigwedge^rW\i \bigwedge^rR((z))^r=R((z))$ is trivial, {\it i.e.}
equal to $R[[z]]\i R((z))$; this is equivalent to saying that
the projective R-module $W/z^N\,R[[z]]^r$  (or $z^{-N}\,R[[z]]^r/W$) is of
rank $Nr$. Then: \medskip

{\pc PROPOSITION} 2.3.-- {\it The $k$-space ${\cal
Q}$ {\rm (resp.} ${\bf GL}_r(K)/{\bf GL}_r({\cal O})${\rm )} represents  the
functor which associates to  a $k$-algebra $R$ the set of special lattices
(resp.  of lattices) $W\i R((z))^r$. The group ${\bf SL}_r(K)$ acts on ${\cal
Q}$ by $(\gamma,W)\mapsto \gamma W$ (for} $\gamma\in GL_r\bigl(R((z))\bigr)$,
$W\i R((z))^r$).

 \ind Let us fix the $k$-algebra $R$, and consider our diagram (1.3)
\vskip -20pt
$$\diagram{
D^*& \lhook\joinrel\mathrel{\hfl{}{}}&D\cr
\vfl{}{}&&\vfl{}{}\cr
X^*&\lhook\joinrel\mathrel{\hfl{}{}}&X\cr
}$$\vskip -10pt
where for simplicity we have dropped the suffix $R$.
Let us start with a pair $(E,\rho)$ over $X$. The trivialization
$\rho$ gives an isomorphism
$R((z))^r\longrightarrow  H^0(D^*,E_{|D^*})$;  the inverse image $W$ of
 $H^0(D,E_{|D})$ is a lattice in $R((z))^r$, and it is a special lattice if
$\wedge^r\rho$ extends to a trivialization of $\bigwedge^rE$ over $X$.
\ind Conversely, given a lattice $W$ in  $R((z))^r$, we define a vector bundle
$E_W$ on $X$ by glueing the trivial bundle over $X^*$ with the bundle
 on $D$ associated to the $R[[z]]$-module $W$; the glueing isomorphism is the
map $W\otimes_{R[[z]]}R((z))\longrightarrow R((z))^r $ induced by the
embedding $W\mono R((z))^r$. By definition $E_W$ has a natural trivialization
$\rho^{}_W$ over $X^*$, and if $W$ is a special lattice $\wedge^r \rho$
extends to a trivialization of $\bigwedge^rE$ over $X$. It is easy to check
that these two constructions are inverse one of each other.

\ind Let $\gamma$ be an element of $GL_r\bigl(R((z))\bigr)$, corresponding
to a triple $(E,\rho,\sigma)$ (1.4). By construction the corresponding
lattice is $\rho^{-1}\sigma\bigl(R[[z]]^r\bigr)=\gamma\bigl(R[[z]]^r\bigr)$.
This proves the last assertion of the proposition. \cqfd
\medskip
\ind Recall that we have denoted by $S^{(N)}$ the subscheme of ${\bf
SL}_r(K)$ parametrizing matrices $A(z)$ such that $A(z)$ and $A(z)^{-1}$ have
a pole of order $\le N$; it is stable under right multiplication by
$S^{(0)}={\bf SL}_r({\cal O})$. We will denote by ${\cal Q}^{(N)}$ its image
in ${\cal Q}$,  i.e. the quotient $k$-space $S^{(N)}/S^{(0)}$. \smallskip

 {\pc COROLLARY} 2.4.--
{\it Let $F_N$ be a free module of rank $r$ over the ring $k[z]/(z^{2N})$
{\rm (}so that $F_N$ is a $k$-vector space of dimension $2rN${\rm ).} The
$k$-space ${\cal Q}^{(N)}=S^{(N)}/S^{(0)}$ is isomorphic to the (projective)
variety  of $rN$-dimensional subspaces $G$ of $F_N$ such that $zG\i G$.}

\ind Let $R$ be a \al a. An element $\gamma$  of $SL_r\bigl(R((z))\bigr)$
belongs to  $S^{(N)}$ if and only if the lattice $W=\gamma R[[z]]^r$ satisfies
$z^N\,R[[z]]^r\i W \i z^{-N}\,R[[z]]^r$. Therefore ${\cal Q}^{(N)}(R)$ is
the subset of ${\cal Q}(R)$ consisting of special lattices with the above
property. Associating to such a lattice its image $G$ in $R\otimes_k
F_N=z^{-N}\,R[[z]]^r/ z^{N}\,R[[z]]^r$ gives a functorial bijection of ${\cal
Q}^{(N)}(R)$ onto the set of direct sub-R-modules $G$ of rank $rN$ of
$R\otimes_k F_N$ which are stable under multiplication by $z$, i.e. an
isomorphism of ${\cal Q}^{(N)}$ onto the subvariety of the Grassmannian ${\bf
G}(rN,F_N)$ corresponding to $z$-stable subspaces of $F_N$. \cqfd \bigskip
\ind We should point out that {\it we do not know whether the varieties
${\cal Q}^{(N)}$ are reduced;} this will cause us some trouble in the sequel.
\smallskip

 \ind Recall that we have denoted by ${\bf SL}_r({\cal O}_-)$ the
subgroup of ${\bf SL}_r(k[z^{-1}])$ parametrizing matrices  $\sum_{n\ge 0}
A_nz^{-n}$ with $A_0=I$. It is an ind-variety. \medskip

{\pc THEOREM} 2.5.-- {\it The $k$-space ${\cal Q}=  {\bf SL}_r(K)/{\bf
SL}_r({\cal O})$ is an ind-variety, direct limit of the system of projective
varieties $({\cal Q}^{(N)})_{N\ge 0}$.  It is covered by open subsets which
are isomorphic to ${\bf SL}_r({\cal O}_-)$, and over which the fibration
$p:{\bf SL}_r(K)\longrightarrow {\cal Q}$ is trivial.}

\ind The first assertion follows from Cor. 2.4. Recall from prop. 1.10 that
the ind-\nobreak group ${\bf SL}_r(K)$  contains an open subset ${\bf
SL}_r(K)^0$ which is isomorphic to\break  ${\bf SL}_r({\cal O}_-)\times  {\bf
SL}_r({\cal O})$,  the isomorphism being equivariant with respect to the
right action of ${\bf SL}_r({\cal O})$. Since ${\bf
SL}_r(K)$ is covered by the open subsets $g{\bf SL}_r(K)^0$ for $g\in
SL_r(K)$, the second assertion follows. \cqfd

\bigskip
\ind The (left) action of  ${\bf SL}_r(K)$ on ${\cal Q}$ restricts to an
action  on the variety ${\cal Q}^{(N)}$ of the group scheme
${\bf SL}_r({\cal O})$ (which actually acts through its finite dimensional
quotient ${\bf SL}_r \bigl(k[z]/(z^{2N})\bigr)$). We are going to study the
orbits of this action. Let us   denote by $\omega$ the class of $I$ in ${\cal
Q}(k)$, and by  $z^{\bf d}$ the matrix ${\rm diag}\,(z^{d_1},\ldots,
z^{d_r})$.
  \medskip
{\pc PROPOSITION} 2.6.-- a)  {\it The orbits of $SL_r({\cal O})$ in ${\cal
Q}(k)$ are the orbits of the points $z^{\bf d}\omega$, where ${\bf d}$ runs
through the sequences $d_1\le \ldots \le d_r$ with $\sum d_i=0$.}
\ind b) {\it The orbit of $z^{\bf d'}\omega$ lies in the closure of the
orbit of $z^{\bf d}\omega$ if and only if one has $d'_1+\ldots+d'_p\ge
d_1+\ldots+d_p$ for $1\le p\le r$.}

\ind c) {\it The subset ${\cal Q}^{(N)}(k)$ is the union of the orbits of
the points $z^{\bf d}\omega$, where ${\bf d}$ runs through the sequences with
$-N\le d_1\le \ldots \le d_r\le N$ and $\sum d_i=0$.}
\ind d) {\it Let ${\bf d}(N)$ denote the sequence $d_1\le \ldots \le d_r$
with $d_i=-N$ for $i<{r+1\over 2}$, $d_i=N$ for $i>{r+1\over
2}$, and $d_{r+1\over 2}=0$ when $r$ is odd. Then the orbit of $z^{{\bf
d}{(N)}}\omega$ is dense in ${\cal Q}^{(N)}$.}

\ind e) {\it The variety ${\cal Q}^{(N)}$ is irreducible.}

\ind Let $W$ be a lattice in $k((z))^r$. Since $k[[z]]$ is a principal ring,
there exists a uniquely determined sequence of integers $d_1\le \ldots \le
d_r$ and a basis $(e_1,\ldots, e_r)$ of  $k[[z]]^r$ such that $W$ is the
lattice spanned by $z^{d_1}e_1,\ldots,z^{d_r}e_r$; this lattice is special if
and only if $\sum d_i =0$. This means that  the point $W$ of ${\cal Q}(k)$
belongs to the orbit of $z^{{\bf d}}\omega$, which proves a). Since the
condition $W\in {\cal Q}^{(N)}$ is equivalent to  $-N\le d_1\le \ldots \le
d_r\le N$, c) follows.

\ind The formula
$$\pmatrix{t^{-1}z&t^{-1}\cr -t& 0}
\pmatrix{
z^{d_1}&0\cr
0& z^{d_2}}
\pmatrix{
t&-t^{-1}z^{d_2-d_1-1}\cr
0&t^{-1}} =\pmatrix{z^{d_1+1}&0\cr
-t^2z^{d_2}& z^{d_2-1}}$$
 for $t\in k^*$, shows that the point $z^{\bf d'}\omega$ belongs to the
closure of the orbit of $z^{\bf d}\omega$ whenever ${\bf d'}$ is obtained
from ${\bf d}$ by replacing a pair of indices $(d_i,d_j)$ with $d_i<d_j$ by
the pair $(d_i+1, d_j-1)$. An easy combinatorial argument  then shows that
every sequence ${\bf d'}$ with $d'_1+\ldots+d'_p\ge d_1+\ldots+d_p$ for $1\le
p\le r$ can be obtained from ${\bf d}$ by iterating this
 operation, which proves the ``if" part of b). To prove the converse, observe
that any matrix  $A(z)$ in $SL_r({\cal O})\,z^{\bf d}\,SL_r({\cal O})$ has
the  property that $z^{d_1}\ldots z^{d_p}$
 divides the coefficients of $\wedge^p A(z)$ for $1\le p\le r$. So if $z^{\bf
d'}$ is a specialization of such a matrix one must have $d_1+\ldots+d_p \le
d'_1+\ldots+d'_p$ for all $p$, which gives b).

\ind The assertion d) is an easy consequence of b) and c). To prove e), it
remains to show that the group scheme  ${\bf SL}_r({\cal O})$ is irreducible;
one way is to observe that the group scheme  ${\bf GL}_r({\cal O})$ is
irreducible (it is isomorphic to $\displaystyle GL_r(k)\times \prod_1^\infty
M_r(k)$, see (1.2)), and maps onto ${\bf SL}_r({\cal O})$ by the morphism
$A\mapsto A\,\delta(A)^{-1}$, where $\delta(A)$ is the diagonal matrix ${\rm
diag}\,(\det A,1,\ldots,1)$. \cqfd \bigskip 

{\it Remark} 2.7.-- One can refine the above decomposition of ${\cal Q}$ as
follows. Let $U$ be the subgroup of $SL_r({\cal O})$ consisting of matrices
$A(z)$ such that $A(0)$ is upper-triangular with diagonal coefficients equal
to $1$. Using the Bruhat decomposition of $SL_r(k)$ one sees easily that the
$SL_r({\cal O})$-orbit of $z^{\bf d}\omega$ is the disjoint union of the sets
$Uz^{{\bf d}_\sigma}\omega$, where ${\bf d}_\sigma$ runs over all
permutations of the sequence $(d_1,\ldots,d_r)$. This is the {\it parabolic
Bruhat decomposition} of the Kac-Moody groups theory [Ku, Sl].

\vskip 1.7cm
{\bf 3. The stack ${\bf SL}_r({\cal O})\backslash{\bf SL}_r(K)/{\bf
 SL}_r(A_X)$}
\smallskip

{\it Stacks}
\ind We will need a few properties of stacks. Rather than giving formal
definitions (for which we refer to [D-M] and especially [L-MB]),  we will try
here to give a rough idea of what stacks are and what they are good for.
Many geometric objects (like vector bundles on a fixed variety, or varieties
of a given type) have no fine moduli space because of  automorphisms. The
remedy is to consider, instead of the set of isomorphism classes, the {\it
groupoid} of such objects (recall that a  groupoid is a category where every
arrow is an isomorphism).

\ind A {\it stack} over $k$  associates to any $k$-algebra $R$ a groupoid
$F(R)$, and to any homomorphism $u:R\rightarrow S$ a functor
$F(u):F(R)\rightarrow F(S)$; these data should satisfy some natural
compatibility conditions as well as some localization properties.

\ind By considering a set as a groupoid (with the identity of
each object as only arrows), a $k$-space can be viewed in a natural way as a
stack over $k$. Conversely, a stack over $k$ with the property that any
object has the identity as only automorphism is a $k$-space.
\bigskip

{\it Examples}.-- (3.1) The {\it moduli stack} ${\cal GL}_X(r)$ of rank $r$
vector bundles on $X$ is defined by associating to a \al a $R$ the groupoid
of rank $r$ vector bundles over $X_R$. Similarly, one defines a  stack ${\cal
SL}_X(r)$ by associating to $R$ the groupoid of pairs $(E,\delta)$, where $E$
is a vector bundle over $X_R$ and $\delta:{\cal O}_{X_R}\rightarrow
\bigwedge^r E$ an isomorphism; this is the fibre over the trivial bundle of
the morphism of stacks $\det:{\cal GL}_X(r)\rightarrow {\cal GL}_X(1)$.

\ind (3.2) Let $Q$ be a $k$-space (1.1), and $\Gamma$ a $k$-group acting on
$Q$. The quotient stack $\Gamma\bk Q$ is defined in the following way: an
object of $F(R)$ is a $\Gamma$-torsor\footnote{\parindent
0.5cm$^1$}{\eightpoint that is a $k$-space over $\Sp(R)$ with an action of
$\Gamma_R$, which after a faithfully flat extension $R\rightarrow R'$ becomes
isomorphic to $\Gamma_{R'}$ acting on itself by multiplication.}\parindent 0cm
together with a $\Gamma$-equivariant morphism $\alpha:P\rightarrow Q$; arrows
in $F(R)$ are defined in the obvious way, and so are  the functors $F(u)$.
The  stack $\Gamma\bk Q$ is indeed the quotient of $Q$ by $\Gamma$ in the
category of stacks,  in the sense that any $\Gamma$-invariant morphism from
$Q$ into a stack $F$ factors through $\Gamma\bk Q$ in a unique way. If
$\Gamma$ acts freely on $Q$ (i.e. $\Gamma(R)$ acts freely on $Q(R)$ for each
\al a $R$), then the stack $\Gamma\bk Q$ is a $k$-space.

\ind When $Q=\Sp(k)$ (with the trivial action), $\Gamma\bk Q$ is the {\it
classifying stack}  $B\Gamma$: for each \al a $R$, $B\Gamma(R)$ is the
groupoid of $\Gamma$-torsors over $\Sp(R)$. \bigskip

{\pc PROPOSITION} 3.4.-- {\it The quotient stack ${\bf
SL}_r(A_X)\backslash{\bf SL}_r(K)/{\bf SL}_r({\cal O})$ is canonically
isomorphic to the algebraic stack ${\cal SL}_X(r)$ of vector bundles on $X$
with trivial determinant. The projection map $\pi: {\bf SL}_r(K)/{\bf
SL}_r({\cal O})\longrightarrow {\cal SL}_r(X)$ is locally trivial for the
Zariski topology. }
 \smallskip
\ind Let us denote as before by ${\cal Q}$ the ind-variety ${\bf
SL}_r(K)/{\bf SL}_r({\cal O})$, and by $\Gamma$ the group ${\bf SL}_r(A_X)$.
The universal vector bundle ${\cal E}$ over $X\times {\cal Q}$, together with
the  trivialization of $\bigwedge^r{\cal E}$ given by $\rho$ (prop. 2.1),
gives rise to a map  $\pi:{\cal Q}\rightarrow {\cal SL}_X(r)$. This map is
$\Gamma$-invariant, hence induces a morphism of stacks
$\overline{\pi}:\Gamma\bk{\cal Q}\rightarrow {\cal SL}_X(r)$.

\ind On the other hand we can define a map ${\cal SL}_X(r)\rightarrow
\Gamma\bk{\cal Q}$ as follows. Let $R$ be a $k$-algebra, $E$ a vector bundle
over $X_R$ and $\delta$ a trivialization of $\bigwedge^rE$. For any
$R$-algebra $S$, let $P(S)$ be the set of trivializations $\rho$ of $E_S$
over $X^*_{S}$ such that $\wedge^r\rho$ coincides with the pull back of
$\delta$.  This defines a $R$-space $P$ on which the group $\Gamma$ acts; by
the lemma below, it is a torsor under $\Gamma$ (it is in fact an ind-scheme,
but we need not  worry about that). To any element of $P(S)$ corresponds a
pair $(E_S,\rho)$, hence by prop. 2.1 an element of ${\cal Q}(S)$. In this
way we associate functorially to an object $(E,\delta)$ of ${\cal
SL}_X(r)(R)$ a $\Gamma$-equivariant map $\alpha:P\rightarrow {\cal Q}$. This
defines a morphism of stacks ${\cal SL}_X(r)\rightarrow \Gamma\bk{\cal Q}$
which is the inverse of $\overline{\pi}$.

\ind The second assertion means that for any scheme $T$ and morphism
\break$f:T\rightarrow {\cal SL}_X(r)$, the pull back to $T$  of the fibration
$\pi$ is (Zariski) locally trivial, i.e. admits local sections. Now $f$
corresponds to a pair $(E,\delta)$, where $E$ is a vector bundle over
$X\times T$ and $\delta$ a trivialization of $\bigwedge^r E$. Let $t\in T$.
By the lemma below, we can find an open neighborhood $U$ of $t$ in  $T$ and a
trivialization $\rho$ of $E_{|X^*\times U}$; modifying $\rho$ by an
automorphism of ${\cal O}^r_{X^*\times U}$ if necessary,  we can moreover
assume $\wedge^r\rho=\delta_{|X^*\times U}$. The pair $(E,\rho)$ defines  a
morphism $g:U\rightarrow {\cal Q}$ (prop. 2.1) such that $\pi\rond g=f$, that
is a section over $U$ of the pull back of the fibration $\pi$. \cqfd
  \bigskip
{\it Lemma} 3.5.-- {\it Let $T$ be a scheme, and $E$ a \vb over $X\times T$
with trivial determinant. Then there exists an open covering $(U_\alpha)$ of
$T$ such that the restriction of $E$ to $X^*\times U_\alpha$ is trivial.}
\smallskip
\ind   We proceed by induction on the rank $r$ of $E$ -- the
case $r=1$ being trivial.  Suppose $r\ge 2$. Let us denote simply by $p$ the
divisor $\{p\}\times T$ in $X\times T$. There exists an  integer $n$ such
that $E(np)$ is spanned by its global sections and has no $H^1$.
 Let $t$ be a point of $T$. An easy count of constants provides a section
$s$ of  $E(np)_{|X\times \{t\}}$ which does not vanish at any point of $X$.
Shrinking $T$ if necessary, we may assume that $s$ is the restriction of a
global section of $E(np)$ which vanishes nowhere on $X\times T$. By
restriction to $X^*\times T$ we get an exact sequence $$0\rightarrow {\cal
O}_{X^*\times T}\longrightarrow E_{|X^*\times T}\longrightarrow F\rightarrow
0\ ,$$ where $F$ is a \vb of rank $r-1$ over $X^*\times T$. Again by
shrinking $T$ if necessary, we may assume that this sequence is split and
(thanks to the  induction hypothesis)  that $F$ is trivial, so $E$ is trivial
over $X^*\times T$. \cqfd

\bigskip
{\it Remark} 3.6.-- The proof of the proposition applies without any
modification to the case of vector bundles with determinant ${\cal O}_X(dp)$,
$d\in {\bf Z}$: the ind-group ${\bf SL}_r(A_X)$ acts by left multiplication
on ${\cal Q}_d\cong {\bf GL}_r(K)^{(d)}/{\bf GL}_r({\cal O})$ (2.2), and the
quotient stack ${\bf SL}_r(A_X)\bk {\cal Q}_d$ is canonically isomorphic to
the moduli stack ${\cal SL}_X(r,d)$ parametrizing vector bundles $E$ on $X_R$
together with an isomorphism ${\cal O}_{X_R}(dp)\iso \bigwedge^rE$. Since
left multiplication by an element $\gamma_d$ of ${\bf GL}_r(K)^{(d)}$ induces
an isomorphism of ${\cal Q}$ onto ${\cal Q}_d$ (2.2), we can also describe
${\cal SL}_X(r,d)$ as the quotient stack $\bigl(\gamma_d^{-1}{\bf
SL}_r(A_X)\gamma_d\bigr)\bk {\cal Q}$. \bigskip

{\it Line bundles over $k$-spaces and stacks}
\ind (3.7)  Let $Q$ be  a $k$-space.
A {\it line bundle} (or a vector bundle, or a coherent
sheaf) ${\cal L}$ on  $Q$ can be defined as the data of a line bundle
(resp. a vector bundle, resp. a coherent sheaf) ${\cal L}_\mu$ on $T$ for
each morphism $\mu$ of a scheme $T$ into $Q$, and of isomorphisms
$g_{\mu,f}:f^*{\cal L}_\mu\iso {\cal L}_{\mu\rond f}$ for each morphism of
schemes $f:T'\rightarrow T$; these data must  satisfy the obvious
compatibility conditions. Morphisms of line bundles (resp. ...)
 are defined in an analogous way; in particular, a section of ${\cal L}$  is
a  compatible family of sections $s_\mu\in H^0(T, {\cal L}_\mu)$,  which means
$g_{\mu,f}(f^*s_\mu)=s_{\mu\rond f}$ for all $f:T'\rightarrow T$. We leave
to the reader to check that all the standard constructions for line bundles
on schemes extend naturally to this situation.

\ind Of course these definitions coincide with the usual ones when $Q$ is a
scheme. Suppose $Q$ is an ind-scheme, limit of an increasing sequence $Q_n$
of schemes; then a line bundle ${\cal L}$ on $Q$ is determined by the data of
a line bundle $(L_n)$ on $Q_n$ for each $n$, and isomorphisms $L_q|_{Q_p}\iso
L_p $ for $q\ge p$, again with  the obvious compatibility conditions. The
space $H^0(Q,{\cal L})$ is then the inverse limit of the system
$\bigl(H^0(Q_n,L_n)\bigr)_{n\ge 1}$.

\ind These definitions can be easily generalized to the case of stacks. A
line bundle ${\cal L}$ on a stack ${\cal S}$  is defined as the data of a
line bundle ${\cal L}_\mu$ on  $T$ for each scheme $T$ and object $\mu$ of
the groupoid ${\cal S}(T)$, and of an isomorphism
$g_\alpha:f^*L_\mu\rightarrow L_\nu$ for each morphism $f:T'\rightarrow T$
and each arrow $\alpha:f^*\mu\rightarrow \nu$ in ${\cal S}(T')$ -- these data
should satisfy some standard compatibility conditions. A section of ${\cal
L}$ is again given by a family of sections $s_\mu\in H^0(T, {\cal L}_\mu)$
such that $g_\alpha(f^*s_\mu)=s_\nu$ for each arrow $\alpha:f^*\mu\rightarrow
\nu$ in ${\cal S}(T')$.

 \bigskip

(3.8) {\it Example: the determinant bundle}.-- Let $T$ be a scheme, and $E$
a vector bundle on $X\times T$. The derived direct image $R(pr^{}_T)_*(E)$ is
given by a complex of vector bundles $L^0\rightarrow  L^1$. The line bundle
$\det (L^1)\otimes \det(L^0)^{-1}$ is independent of the choice of this
complex, hence canonically defined on $T$; this is the ``determinant of the
cohomology" $\det R\Gamma_T(E)$.
 Associating to each  bundle $E$ on $X\times T$ the
line bundle $\det R\Gamma_T(E)$ defines a line bundle ${\cal L}$ on the stack
${\cal SL}_X(r)$ (or ${\cal GL}_X(r)$), the {\it determinant line bundle.}
\ind There is a useful way to produce sections of the line bundle $\det
R\Gamma_T(E)$ and of its multiples.  Suppose for simplicity that $T$ is
integral, and that the line bundle $\bigwedge^rE$ is the pull back of some
line bundle on $X$. Let $F$ be a vector bundle of rank $s$ on $X$; let us use
the same notation to denote  its pull back to $X\times T$. Then {\it the line
bundle $\det R\Gamma_T(E\otimes F)$ is isomorphic to $\bigl(\det
R\Gamma_T(E)\bigr)^{\otimes s}$.} (write $F$ as an extension to reduce to the
case $s=1$, then use repeteadly the exact sequence $0\rightarrow E\rightarrow
E(q)\rightarrow E\otimes ({\cal O}(q)/{\cal O})\rightarrow 0$ to prove that
the line bundle $\det R\Gamma_T(E(D))$ is isomorphic to $\det R\Gamma_T(E)$
for any divisor $D$ on $X$). Put $E_t:=E_{|X\times\{t\}}$ for $t\in T$.
Choose $F$ such that the vector bundle $E_t\otimes F$ has trivial cohomology
for some $t$ in $T$. Let $L^0\phfl{u}{}L^1$ be a complex of vector bundles
isomorphic to $R(pr^{}_T)_*(E\otimes F)$. Then  $\det u$ is a nonzero section
of  $\det R\Gamma_T(E\otimes F)\cong \bigl(\det R\Gamma_T(E)\bigr)^{\otimes
s}$, which is well defined up to an invertible function on $T$. In
particular, its divisor $\Theta_F$ is canonically defined on $T$ (the support
of $\Theta_F$ is the set of points $t\in T$ such that $H^0(X,E_t\otimes
F)\not=0$).

 \bigskip

{\it Example} 3.9.-- Let $G$ be a $k$-group,  $H$ a $k$-subgroup of $G$, and
$\chi:H\rightarrow {\bf G}_m$ a character of $H$. As usual we  associate to
this situation  a line bundle ${\cal L}_\chi$ on $G/H$: the group $H$ acts
freely on the trivial bundle $G\times {\bf A}^1$ by $h\cdot (g,t)=(gh\,
,\,\chi(h^{-1})t)$ and we define ${\cal L}_\chi$ as the quotient $k$-space
$(G\times  {\bf A}^1)/H$. It is an easy exercise in descent theory to prove
that the pull back of the fibration ${\cal L}_\chi\rightarrow G/H$ by any
morphism $\mu:T\rightarrow G/H$ is indeed a line bundle (use the fact that
$\mu$ lifts locally to $G$, and that the pull back of ${\cal L}_\chi$ to $G$
is the trivial bundle). Again by descent, sections of ${\cal L}_\chi$
correspond in a one-to-one way to sections of the trivial bundle $G\times
{\bf A}^1$ over $G$ which are $H$-invariant, that is to {\it functions $f$ on
$G$ such that $f(gh)=\chi(h^{-1})f(g)$} for any \al a $R$ and elements $g\in
G(R)$, $h\in H(R)$.

  \ind There is a more fancy way to describe the line bundle ${\cal L}_\chi$.
Consider the classifying stack $B{\bf G}_m$ over $k$ (3.2). A morphism
$\mu:T\rightarrow B{\bf G}_m$ is given by a ${\bf G}_m$-torsor over $T$,
which defines a line bundle ${\cal U}_\mu$ over $T$: this defines the
universal line bundle ${\cal U}$ over $B{\bf G}_m$. As a stack over $B{\bf
G}_m$, it is simply the quotient ${\bf A}^1/{\bf G}_m$. The character
$\chi:H\rightarrow {\bf G}_m$ induces a morphism $B\chi:BH\rightarrow B{\bf
G}_m$, hence a line bundle $(B\chi)^*{\cal U}$ on $BH$. The structural map
$G\rightarrow \Sp(k)$ induces a morphism of stacks $G/H\rightarrow BH$, hence
by pull back we get a line bundle on $G/H$, which is (almost by definition)
${\cal L}_\chi$.

\ind Let $G'$ be another $k$-group,  $H'$ a $k$-subgroup of $G'$, and
$f:G'\rightarrow G$ a morphism of $k$-groups which maps $H'$ into $H$.  It
follows from either of these definitions that the pull back of ${\cal
L}_\chi$ by the morphism $G'/H'\rightarrow G/H$ induced by $f$  is the line
bundle ${\cal L}_{\chi'}$ associated to the character $\chi':=\chi\rond f$ of
$H'$.

\vskip 1.7cm
{\bf 4. The central extension}
\smallskip
\ind Let $\pi:{\cal Q}\longrightarrow {\cal SL}_X(r)$ be the canonical
morphism of stacks defined in the preceding section. We want to identify the
line bundle $\pi^*{\cal L}$ on  ${\cal Q}$.

\ind It will turn out that, though it is invariant under the action of ${\bf
SL}_r(K)$, this line bundle does {\it not} admit an action of ${\bf
SL}_r(K)$. But it does admit an action of a canonical extension
$\widehat{\bf SL}_r(K)$ of ${\bf SL}_r(K)$, which we are now going to
describe paraphrasing [S-W].
 \medskip

{\it The canonical extension of the Fredholm group}
\smallskip
\ind (4.1) Let $V$ be an infinite-dimensional vector space over $k$. Denote
by $\End^f(V)$ the two-sided ideal of $\End(V)$ formed by the endomorphisms
of finite rank and  by
 ${\cal F}(V)$ the group of units of the
quotient algebra $\End(V)/\End^f(V)$.  The elements of ${\cal F}(V)$ are
classes of equivalence of endomorphisms with finite-dimensional  kernel and
cokernel. We let ${\cal F}(V)^0$ be the subgroup of classes of index $0$
endomorphisms, {\it i.e.} endomorphisms with $\dim\Ker u =\dim\Coker u$.
 It is an easy exercise
to show that the image of the canonical homomorphism $\Aut(V)\rightarrow
{\cal F}(V)$ is ${\cal F}(V)^0$;  its kernel  consists of the automorphisms
$u$ of $V$ such that $u\equiv I$ (mod.~$\End^f(V)$). The determinant of such
an endomorphism is naturally defined, e.g. by the formula $\displaystyle
\det(I+v)=\sum_{n\ge 0}\Tr\wedge^nv$. Let us denote by
$\bigl(I+\End^f(V)\bigr)_1$  the subgroup of automorphisms of the form $I+v$
with $v\in \End^f(V)$ and $\det(I+v)=1$; we get an exact sequence
$$1\rightarrow
k^*\longrightarrow \Aut(V)/\bigl(I+\End^f(V)\bigr)_1\longrightarrow {\cal
F}(V)^0 \rightarrow 1\ .$$  For $v\in\End^f(V),\ u\in\Aut(V)$, one has
$\det(I+uvu^{-1})=\det(I+v)$; this means that the element $I+v$ belongs to
the {\it center} of the group $\Aut(V)/\bigl(I+\End^f(V)\bigr)_1$. We have
thus defined  a canonical central extension of the group ${\cal F}(V)^0$ by
$k^*$.  \medskip

\ind (4.2) We want to make sense of this in an algebraic setting, at least
at the level of
 $k$-groups. We define the $k$-space
$\End(V)$ in an obvious way, as the functor $R\mapsto \End_R(V\otimes_k
R)$, and the $k$-group $\Aut(V)$ as its group of units. We'll say that an
endomorphism of $V\otimes R$ has finite rank if  its image is contained in a
finitely generated submodule; we define $\End^f(V)(R)$ as the ideal formed by
these endomorphisms, and take for ${\cal F}(V)$ the group of units of the
algebra $\End(V)/\End^f(V)$. We don't know a good definition for the subgroup
${\cal F}(V)^0$, so we just define it as the image of $\Aut(V)$ in ${\cal
F}(V)$. We then get again a central extension of $k$-groups $$1\rightarrow
{\bf G}_m\longrightarrow \Aut(V)/\bigl(I+\End^f(V)\bigr)_1\longrightarrow
{\cal F}(V)^0 \rightarrow 1\ .\leqno({\cal F})$$ \bigskip

{\it The central extension of ${\bf SL}_r(K)$}
\smallskip
\ind Let us go back to the ind-group ${\bf GL}_r(K)$. We choose a supplement
$V$ of ${\cal O}^r$ in $K^r$. For any $k$-algebra, we get a direct
sum decomposition (over $k$)
$$R((z))^r=V_R\oplus R[[z]]^r\ ,$$
with $V_R:=V\otimes_k R$.  Let $\gamma$
be an element of  $GL_r\bigl(R((z))\bigr)$, and let $$\gamma =
\pmatrix{a(\gamma)& b(\gamma)\cr c(\gamma)& d(\gamma)\cr}\leqno(4.2)$$ be its
matrix with respect to the above decomposition.  Let $\overline{a}(\gamma)$
denote the class of $a(\gamma)$ in $\End(V_R)/\End^f(V_R)$.

\medskip

{\pc PROPOSITION} 4.3.--  a) {\it The map $\gamma\mapsto \overline{a}
(\gamma)$ is a group homomorphism from $GL_r\bigl(R((z))\bigr)$ into} ${\cal
F}(V_R)$; {\it it defines a morphism of $k$-groups:} \vskip -10pt
$$\overline{a}:{\bf GL}_r(K)\rightarrow {\cal F}(V)\ .$$

\ind  b) {\it Let $V'$ be another supplement of ${\cal O}^r$ in $K^r$, and
let $\overline{a'} :{\bf GL}_r(K)\rightarrow {\cal F}(V')$  be the morphism
associated to $V'$. Let $\varphi:V\rightarrow V'$ be the isomorphism obtained
by restricting to $V$  the projector onto $V'$. Then   $\overline{a'}$ is
equal to $\varphi\,\overline{a}\,\varphi^{-1}$.} \ind   Since $\gamma$ maps
$R[[z]]^r$ into $z^{-N}R[[z]]^r$ for some $N$, the map
$b(\gamma):R[[z]]^r\rightarrow V_R$ has finite rank. From this and the
formula for the product of two matrices follows first that the endomorphism
$a$ of $V_R$ is invertible modulo finite rank endomorphisms,  then that the
map  $GL\bigl(R((z))\bigr)\longrightarrow {\cal F}(V)(R)$ which associates to
$\gamma$ the class of $a(\gamma)$ is a group homomorphism. This proves a).
\smallskip

\ind Let $p$, $q$ be the projectors of $R((z))^r$ onto $V'_R$ and $R[[z]]^r$
relative to the decomposition  $R((z))^r=V'_R\oplus R[[z]]^r$, and let
$\pmatrix{a'(\gamma)& b'(\gamma)\cr c'(\gamma)& d'(\gamma)\cr}$ be the matrix
of an element $\gamma\in GL_r\bigl(R((z))\bigr)$ relative to this
decomposition. An easy computation gives
$p\,a(\gamma)=a'(\gamma)\,p+b'(\gamma)\,q$. Since $b'(\gamma)$ has finite
rank, we get the equality $a'(\gamma)\equiv\varphi a(\gamma)\varphi^{-1}$
(mod. ${\End^f(V_R)}$). \cqfd

 \medskip
{\pc PROPOSITION} 4.4.-- {\it  Let $R$ be a
\al a and $\gamma$ an element of $SL_r\bigl(R((z))\bigr)${\rm ;} locally on
$\Sp(R)$ (for the Zariski topology), the endomorphism $a(\gamma)$
of $V_R$ is equivalent mod. $\End^f(V_R)$ to an automorphism.}  \smallskip

\ind By prop. 4.3 b), it is enough to prove the result for one particular
choice of $V$; we'll take
 $V=\bigl(z^{-1}k[z^{-1}]\bigr)^r$. The assertion
 is clear when $\gamma$ belongs to $SL_r\bigl(R[[z]]\bigr)$ or to
$SL_r\bigl(R[z^{-1}]\bigr)$: in those cases the matrix (4.2) is  triangular,
so that $a(\gamma)$ itself is an isomorphism. The result then follows when
$R$ is a field, since any matrix $\gamma\in SL_r\bigl(R((z))\bigr)$ can be
written as  a product of elementary matrices $I+\lambda E_{ij}$, where
$\lambda$ can be taken either in $R[[z]]$ or in $R[z^{-1}]$. The general case
is a consequence of
 the following lemma:
\medskip
{\it Lemma} 4.5.-- {\it Locally over $\Sp(R)$, any element  $\gamma$ of
$SL_r\bigl(R((z))\bigr)$ can be written $\gamma_0\gamma^-\gamma^+$, with
$\gamma_0\in SL_r(K)$, $\gamma^-\in SL_r(R[z^{-1}])$, $\gamma^+\in
SL_r(R[[z]])$.} \ind Let us assume first that the \al a $R$ is finitely
generated. Let $t$ be a closed point of $\Sp(R)$; put $\gamma_0=\gamma(t)$.
By (1.12) $\gamma_0^{-1}\gamma$ can be written in a neighborhood of $t$ as
$\gamma^-\gamma^+$, hence the result in this case. \ind In the general case,
$R$ is the  union of its finitely generated subalgebras $R_\alpha$. Let
$p:{\bf SL}_r(K)\longrightarrow {\cal Q}={\bf SL}_r(K)/{\bf SL}_r({\cal O})$
be the quotient map. Since ${\cal Q}$ is an ind-variety, the morphism $p\rond
\gamma:\Sp(R)\rightarrow {\cal Q}$  factors through $\Sp(R_\alpha)$ for some
$\alpha$. Locally over $\Sp(R_\alpha)$, this morphism can be written $p\rond
\gamma_\alpha$ for some element $\gamma_\alpha$ of
$SL_r\bigl(R_\alpha((z))\bigr)$, which differ from $\gamma$ by an element of
$SL_r(R[[z]])$ (thm. 2.5). Since $R_\alpha$ is of finite type, the lemma
holds for $\gamma_\alpha$, hence also for $\gamma$. \cqfd \bigskip

{\pc COROLLARY} 4.6.-- {\it The image of ${\bf
SL}_r(K)$ by $\overline{a}$ is contained in the subgroup ${\cal F}(V)^0$.
\cqfd}
\medskip

\ind We will denote by $\widehat{\bf SL}_r(K)\rightarrow {\bf SL}_r(K)$ the
pull back of the central extension $({\cal F})$ by $\overline{a}$, so that we
get a central extension of $k$-groups $$0\longrightarrow  {\bf
G}_m\phfl{}{}\widehat {\bf SL}_r(K)\phfl{\psi}{}{\bf SL}_r(K)\longrightarrow
0 \ .\leqno{({\cal E})}$$
 By descent theory any ${\bf
G}_m$-torsor over a scheme is representable, so the $k$-group $\widehat{\bf
SL}_r(K)$ is also an ind-group.
\smallskip

\ind (4.7) Let $R$ be a $k$-algebra; an element of $\widehat {\bf
SL}_r(K)(R)$ is given, locally on $\Sp R$, by a pair $(\gamma,u)$ with
$\gamma$ in $SL_r(R((z)))$, $u$ in $\Aut(V_R)$, and $u\equiv a(\gamma)$ (mod.
$\End^f(V_R)$); two pairs $(\gamma,u)$ and $(\gamma,v)$ give the same element
if $u^{-1}v$ (which belongs to $I+\End^f(V_R)$) has determinant $1$. In
particular, the kernel of $\psi$ consists of the pairs $(I,u)$ with $u\in
I+\End^f(V)$, modulo the pairs $(I,u)$ with $\det u=1$; the map $u\mapsto
\det u$ provides an isomorphism from  $\Ker \psi$ onto ${\bf G}_m$.

\ind Because of prop. 4.3 b), the extension $({\cal E})$ is independant of
the choice of the supplement $V$ of ${\cal O}^r$ in $K^r$. More precisely,
given two such supplements $V$ and $V'$, there is a canonical isomorphism
from the group $\widehat {\bf SL}_r(K)$ defined using $V$ onto the group
$\widehat {\bf SL}'_r(K)$ defined using $V'$: it associates to a pair
$(\gamma,u)$ as above the pair $(\gamma,\varphi u\varphi^{-1})$, where
$\varphi$ is the natural isomorphism from $V$ onto $V'$ ({\it loc. cit.}).
One can then define in the usual way a {\it canonical} central extension of
${\bf SL}_r(K)$ by taking the projective (or inductive) limit, over the set
of all supplements  of ${\cal O}^r$ in $K^r$, of the extensions we have
constructed.   \medskip

\ind (4.8) Let $H$ be a sub-$k$-group of ${\bf SL}_r(K)$, such that  ${\cal
O}^r$ (resp. $V$) is stable under $H$. Then {\it the extension $({\cal E})$
is canonically split over $H$}. For any element $\gamma$ of $H(R)$ satisfies
$b(\gamma)=0$ (resp. $c(\gamma)=0$), so that the map $\gamma\mapsto
a(\gamma)$ is a  homomorphism from  $H(R)$ into $\Aut(V_R)$. Then the map
$\gamma\mapsto (\gamma,a(\gamma))$ defines a  section of $\psi$ over $H$. In
particular, we see that {\it the pull back $\widehat {\bf SL}_r({\cal O})$
of ${\bf SL}_r({\cal O})$ is canonically isomorphic to ${\bf SL}_r({\cal
O})\times {\bf G}_m$}. We will denote by $\chi_0: \widehat {\bf SL}_r({\cal
O})\longrightarrow {\bf G}_m$ the second projection; if the element
$\tilde\delta$ of $\widehat {\bf SL}_r({\cal O})(R)$ is represented by a pair
$(\delta,v)$, one has  $\chi_0(\tilde \delta)=\det(a(\delta)^{-1}v)$. \ind
More generally, suppose that there exists an element $\lambda\in SL_r(K)$
such that the subgroup $H$ preserves the subspace $\lambda({\cal O}^r)$
(resp. $\lambda(V)$).  We choose an automorphism $u$ of $V$ such that
$u\equiv a(\lambda)$ (mod. $\End^f(V)$), and define a section of $\psi$ over
$H$ by $\gamma\mapsto ua(\lambda^{-1}\gamma\lambda)u^{-1}$. This section is
independent of the choice of $u$, so once again the group $H$ embeds
canonically into $\widehat{\bf SL}_r(K)$.

\bigskip
{\it The Lie algebra of the central extension and the Tate residue}
\smallskip
\ind (4.9) We want to show that at the level of Lie algebras, the extension
$({\cal E})$ is the universal central extension which appears in the theory
of Kac-Moody algebras  [K]. This is essentially known (see e.g. [A-D-K],
where very similar computations appear). We have included the computation
because it is extremely simple  and gives a nice generalization of the
residue defined by Tate in [T].

\ind Let us start from the central extension $({\cal F})$. Since ${\cal
F}(V)$ is the group of invertible elements of the associative algebra
$\End(V)/\End^f(V)$, its
 Lie algebra  is simply the quotient of the Lie algebra $\End(V)$ by the
ideal $\End^f(V)$. The Lie algebra of $(I+\End^f(V))_1$ is the sub-Lie
algebra $\End^f(V)_0$ of $\End^f(V)$ consisting of traceless endomorphisms.
Therefore the Lie algebras extension corresponding to $({\cal F})$ is
$$0\rightarrow k\longrightarrow \End(V)/\End^f(V)_0\longrightarrow
\End(V)/\End^f(V)\rightarrow 0\ .\leqno({\goth F})$$

\ind Let $\alpha$ be an element of ${\goth sl}_r\bigl(k((z))\bigr)$; it
corresponds to the element $I+\varepsilon\alpha$ of
$SL_r\bigl(k[\varepsilon]((z))\bigr)$. Since
$a(I+\varepsilon\alpha)=I+\varepsilon a(\alpha)$, the tangent map
$L(\overline{a})$ at $I$ to $\overline{a}: {\bf SL}_r(K)\rightarrow {\cal
F}(V)$ associates to $\alpha$ the class of $a(\alpha)$ in
$\End(V)/\End^f(V)$. By construction the extension of ${\goth sl}_r(K)$ we
are looking for is the pull back of ${\goth F}$ by $L(\overline{a})$. This
means that the
 Lie algebra  $\widehat {\goth sl}_r(K)$ of
 $\widehat{\bf SL}_r(K)$ consists of pairs $(\alpha,u)$ with
$\alpha\in{\goth sl}_r(K)$, $u\in\End(V)$, $a(\alpha)\equiv u$ (mod.
$\End^f(V)$); two  pairs $(\alpha,u)$ and $(\alpha,v)$ give the same element
if $\Tr(u-v)=0$. We get a central extension $$0\longrightarrow
k\phfl{}{}\widehat {\goth sl}_r(K)\phfl{L(\psi)}{} {\goth
sl}_r(K)\longrightarrow 0 \leqno({\goth E})$$ with
$L(\psi)(\alpha,u)=\alpha$, the kernel of $L(\psi)$ being identified with $k$
by $(0,u)\mapsto\Tr u$. As before, this extension does not depend on the
choice of the supplement $V$ of ${\cal O}^r$ in $K^r$. We claim that it is
the well-known {\it universal central extension} of $ {\goth sl}_r(K)$.
Recall that this extension is obtained by  defining a Lie algebra structure
on ${\goth sl}_r(K)\oplus k$ by the formula
$$\bigl[(\alpha,s),(\beta,t)\bigr]=\bigl([\alpha,\beta]\,,\Res_0\Tr({d\alpha\over
dz}\,\beta)\bigr)\ ;$$ the projection $p$ onto the first summand and the
injection $i$ of the second one  define the universal central extension
$$0\rightarrow k\phfl{i}{} {\goth sl}_r(K)\oplus k \phfl{p}{} {\goth
sl}_r(K)\rightarrow 0\ .\leqno({\goth U})$$  \medskip  {\pc PROPOSITION}
4.10.-- {\it There exists a Lie algebra isomorphism  $\widehat {\goth
sl}_r(K)\iso  {\goth sl}_r(K)\oplus k$ inducing an isomorphism of the
extension ${\goth E}$ onto the universal central extension of ${\goth
sl}_r(K)$.}

\ind It is enough to prove the proposition for one particular choice of $V$;
it will be convenient to choose for $V$ the subspace
${\cal O}_-^r$, where
 ${\cal O}_-$ denotes the
subspace  $z^{-1}k[z^{-1}]$ of $K$.

 \ind Let us define a map
$\varphi:\widehat {\goth sl}_r(K)\longrightarrow {\goth sl}_r(K)\oplus k$
by $\varphi(\alpha,u)= \bigl(\alpha,\Tr(u-a(\alpha))\bigr)$.
 One has $p\rond \varphi=L(\psi)$ and $\varphi$ induces an isomorphism of
$\Ker L(\psi)$ onto $i(k)$; this implies that  $\varphi$ is bijective. It
remains to check  that $\varphi$ is a Lie algebra homomorphism. Since
$\varphi$ maps $\Ker L(\psi)$ into the center $i(k)$, it is enough to prove
the equality $\varphi([\tilde \alpha,\tilde \beta])=[\varphi(\tilde
\alpha),\varphi(\tilde \beta)]$ for $\tilde \alpha=(\alpha,a(\alpha))$,
$\tilde \beta=(\beta,a(\beta))$. This amounts to the following formula:
\medskip

{\it Lemma} 4.11.-- {\it Let $\alpha$, $\beta$ be
two matrices in $M_r(K)$. One has} $$
\Tr\bigl([a(\alpha),a(\beta)]-a([\alpha,\beta])\bigr)=
\Res_0\Tr({d\alpha\over dz}\beta)\ .$$ \smallskip

 \ind This  is precisely Tate's definition of the
 residue [T] in the case $r=1$; we will actually reduce the proof to the
rank $1$ case. \ind   Assume first that for some integer $N$ one has
$\alpha\in z^{N+1}M_r({\cal O})$, \break
 $\beta\in z^{-N}M_r({\cal O})$. For $p\ge 0$, let us denote by $V_p$ the
subspace  $V\cap z^{-p}{\cal O}^r$ of $V$.
 Then $(V_p)_{p\ge 0}$ is
an increasing filtration of $V$, and
 for each $p\ge 0$, the endomorphism
$[a(\alpha),a(\beta)]-a([\alpha,\beta])$ maps $V_p$ into $V_{p-1}$. This
implies that its trace is zero, which gives the formula in this case.

\ind  By bilinearity, we can therefore assume that $\alpha$ and $\beta$ are
polynomial in $z$, and even of the form $z^pA$ for some integer $p$ and some
matrix $A\in M_r(k)$. Let us identify $K^r$ with $K\otimes_k k^r$. The direct
sum decomposition $K^r=V\oplus {\cal O}^r$ is induced by tensor product from
the decomposition $K={\cal O}_-\oplus {\cal O}$. It follows that
$a(z^p\otimes A)$ is $a_1(z^p)\otimes A$, where $a_1(z^p)$ is the
endomorphism of  ${\cal O}_-$ associated to $z^p$. Since the trace of
$u\otimes M$, for $u\in\End^f({\cal O}_-)$ and $M\in M_r(k)$, is $(\Tr
u)\,(\Tr M)$, and since the endomorphisms $z^p$ and $z^q$ of $K$ commute, we
obtain $$\Tr\Bigl([a(z^pA),a(z^qB)]-a([z^pA,z^qB]\Bigr)=\Tr AB\
\Tr\,[a_1(z^p),a_1(z^q)]\ .$$ It remains to compute the trace of the (finite
rank) endomorphism $u=[a_1(z^p),a_1(z^q)]$ of ${\cal O}_-$,
 which we do using the basis $(z^{-n})_{n\ge 1}$ of
${\cal O}_-$. One has $u(z^{-n})=\varepsilon z^{p+q-n}$, with
$\varepsilon\in \{-1,0,1\}$; therefore  $\Tr u$ is zero  except when $p+q=0$.
Assume $q=-p$ and, say, $p\ge 0$; then we find $u(z^{-n})=0$ for $n>p$, and
$u(z^{-n})=z^{-n}$ for $1\le n\le p$. We conclude that $\Tr
\,[a_1(z^p),a_1(z^q)]=\delta_{p,-q}\ p=\Res_0(pz^{p-1}z^q)$, from which the
proposition follows. \cqfd

\bigskip
\ind (4.12)  The above
computations extend in a straightforward way when the base field $k$ is
replaced by an arbitrary \al a $R$. In particular, the kernel of the
homomorphism $\widehat{\bf SL}_r(K)(R[\varepsilon])\longrightarrow
\widehat{\bf SL}_r(K)(R)$ is the
 Lie algebra  $\widehat{\goth sl}_r\bigl(R((z))\bigr)={\goth
sl}_r\bigl(R((z))\bigr)\oplus R$, where the Lie bracket is defined by
formula (4.10).  This defines an adjoint action of the group
$SL_r\bigl(R((z))\bigr)$ onto  $\widehat{\goth sl}_r\bigl(R((z))\bigr)$,
which is trivial on the center  and induces on the quotient ${\goth
sl}_r\bigl(R((z))\bigr)$ the usual  action by conjugation. We claim that it
is given by the following formula:  $$\Ad(\gamma)\
(\alpha,s)=\bigl(\gamma\alpha\gamma^{-1}\,,\,s+\Res_0\Tr(\gamma^{-1}{d\gamma\over
dz}\ \alpha)\bigr)\ .\leqno(4.12)$$ In fact, let  $\gamma\in
SL_r\bigl(R((z))\bigr)$, and let $\ell$ be a $R$-linear form on ${\goth
sl}_r\bigl(R((z))\bigr)$. The condition for the map  $(\alpha,s)\mapsto
(\gamma\alpha\gamma^{-1},s+\ell(\alpha))$ to be a Lie algebra homomorphism is
$\ell([\alpha,\beta])=\displaystyle
\Res_0\Tr\Bigl({d(\gamma\alpha\gamma^{-1})\over
dz}\,\gamma\beta\gamma^{-1}-{d\alpha\over dz}\,\beta\Bigr)$. Since the Lie
algebra  ${\goth sl}_r(S)$, for any ring $S$, is equal to  its commutator
algebra, this condition determines $\ell$ uniquely. On the other hand, it is
checked readily that the linear form
$\alpha\mapsto\Res_0\Tr(\gamma^{-1}{d\gamma\over dz}\ \alpha)$ has the
required property.
  \bigskip

{\it The $\tau$ function}
\smallskip

\ind Let $R$ be a \al a, and $\tilde\gamma$ an element of $\widehat{\bf
SL}_r(K)(R)$. Locally on $\Sp(R)$ we can write $\tilde\gamma=(\gamma,u)$ with
$\gamma$ in $SL_r(R((z)))$, $u$ in $\Aut(V_R)$, and $u\equiv a(\gamma)$ (mod.
$\End^f(V_R)$).  We associate to this pair the element $\tau^{}_V(\gamma,u):=
\det(u\,a(\gamma^{-1}))$ of $R$. This is clearly well-defined, so we get  an
algebraic function $\tau^{}_V$ on $\widehat {\bf SL}_r(K)$.
 \bigskip

{\pc PROPOSITION} 4.13.-- {\it Let $R$ be a $k$-algebra,  $\tilde\gamma$ an
element of  $\widehat{\bf SL}_r(K)(R)$, $\gamma$ its image in
$SL_r\bigl(R((z))\bigr)$.
 One has
$\tau^{}_V(\tilde\gamma\tilde\delta)=\chi_0(\tilde\delta)\tau^{}_V
(\tilde\gamma)$ for all $\tilde\delta$ in $\widehat{\bf SL}_r({\cal O})(R)$.}

\ind Let us choose representatives $(\gamma,u)$ of $\tilde\gamma$ and
$(\delta,v)$ of $\tilde\delta$. Since $b(\delta^{-1})=0$, one has
$a(\delta^{-1}\gamma^{-1})=a(\delta^{-1})a(\gamma^{-1})$, and
\vskip -10pt
$$\tau^{}_V(\tilde\gamma\tilde\delta)=
\det\bigl(uva(\delta^{-1})a(\gamma^{-1})\bigr)=
\det\bigl(va(\delta^{-1})\bigr) \det\bigl(ua(\gamma^{-1})\bigr)=\chi_0(\tilde
\delta)\tau^{}_V(\tilde \gamma)\ .\ \ \carre$$

 \bigskip
\ind (4.14) Let us denote by $\chi$ the character  $\chi_0^{-1}$ of
$\widehat{\bf SL}_r({\cal O})$. The function $\tau^{}_V$  thus defines a
section  of the line bundle ${\cal L}_{\chi}$ on the ind-variety ${\cal
Q}=\widehat{\bf SL}_r(K)/ \widehat{\bf SL}_r({\cal O})$ (3.8). More
generally, let   $\delta\in SL_r(K)$, and let $\tilde \delta$  be a lift of
$\delta$ in $\widehat{SL}_r(K)$; the function $\tilde
\gamma\mapsto\tau^{}_V(\tilde \delta^{-1}\tilde \gamma)$ still defines an
element of $H^0({\cal Q},{\cal L}_\chi)$, whose divisor is
$\delta(\div(\tau^{}_V))$.

\medskip
\ind To conclude this section, let us mention that one gets slightly more
natural conventions by having the group $SL_r({\cal O})$ acting {\it on the
left} on $SL_r(K)$: in particular the twist $\gamma\mapsto \gamma^{-1}$ in
the definition of the $\tau$ function  disappears, and the $\tau$ function
becomes a section of ${\cal L}_{\chi^{}_0}$. We have chosen instead to follow
the standard conventions of Kac-Moody theory.

\vskip 1.7cm
{\bf 5. The determinant bundle}
\smallskip
\ind We will now compare the pull back over
${\cal Q}$ of the determinant line bundle ${\cal L}$ on the moduli stack
with the line bundle ${\cal L}_\chi$ we have just described.
\smallskip
{\pc PROPOSITION} 5.1.-- {\it Let $R$ be a \al a, $\gamma$  an
element of $GL_r(R)$, and  $(E,\rho,\sigma)$  the corresponding triple over
$X_R$ {\rm (1.4)}. There
is a canonical exact sequence
$$0\rightarrow  H^0(X_R,E)\longrightarrow A_X^r\otimes_k
R\phfl{\overline{\gamma}}{} (R((z))/R[[z]])^r \longrightarrow
H^1(X_R,E)\rightarrow 0$$ where $\overline{\gamma}$ is the composition of the
injection $A_X^r\otimes_kR\mono R((z))^r$ deduced from the restriction map
$A_X\rightarrow k((z))$, the automorphism $\gamma^{-1}:R((z))^r\rightarrow
R((z))^r$, and the canonical surjection $R((z))^r\rightarrow
(R((z))/R[[z]])^r$.}
\smallskip
\ind This is in fact the cohomology exact
sequence associated to the exact sequence defined in (1.4) $$0\rightarrow
E\longrightarrow j_*{\cal O}_{X^*}^r\phfl{\overline\gamma}{} f_*({\cal
K}_D/{\cal O}_D)^r\rightarrow 0\ .\quad\carre$$
\ind (5.2) Let us choose an element $\gamma_0$ in $GL_r(K)$ such that the
associated bundle $E_{\gamma_0}$ has trivial cohomology. According to what we
have just seen, this means that the map $\overline{\gamma}_0:A_X^r
\longrightarrow (K/{\cal O})^r$ is an isomorphism, or in other words that
{\it the subspace $V:=\gamma_0^{-1}(A_X^r)$ is a supplement of ${\cal O}^r$
in $K^r$.} Let us identify $A_X^r$ to $V$ with the help of $\gamma_0$, and
the quotient map $K^r\longrightarrow (K/{\cal O})^r$ to the projection of
$K^r$ onto $V$; we obtain that $\overline{\gamma}$ is the composition of the
mappings $V\mono K^r\ \hfl{\gamma^{-1}\gamma_0}{}\ K^r\phfl{}{}V$. In other
words, $\overline{\gamma}$ is the coefficient $a(\gamma^{-1}\gamma_0)$ of the
matrix of $\gamma^{-1}\gamma_0$ with respect to the decomposition
$K^r=V\oplus {\cal O}^r$ (\S 4). We have thus obtained:
\medskip
{\pc PROPOSITION} 5.2.-- {\it Let $\gamma$ be an element of
$GL_r\bigl(R((z))\bigr)$, and let $E$ be the associated vector bundle over
$X_R$. There is a canonical exact sequence}\vskip3pt $$0\rightarrow
H^0(X_R,E)\longrightarrow V_R\ \ghfl{a(\gamma^{-1}\gamma_0)}{}\
V_R\longrightarrow  H^1(X_R,E)\rightarrow 0\ .\quad\carre\leqno(5.2)$$
\bigskip {\pc COROLLARY} 5.3.-- {\it Assume that there exists an automorphism
$u$ of $V_R$ such that} $u\equiv a(\gamma_0^{-1}\gamma)$ (mod.
$\End^f(V_R)$). {\it Then there is an exact sequence}\vskip2pt $$0\rightarrow
H^0(X_R,E)\longrightarrow V_0\hfl{v_0}{}V_0\longrightarrow
H^1(X_R,E)\rightarrow 0\leqno(5.3)$$ {\it where $V_0$ is a free finitely
generated $R$-module, and $\det(v_0)=\tau^{}_V(\gamma_0^{-1}\gamma, u)$.}

 \ind Let $v=u\,a(\gamma^{-1}\gamma_0)\in I+\End^f(V_R)$, and let $V_0$ be a
free finitely generated direct factor of $V_R$ containing $\im (v-I)$. We
denote by  $v_0$  the restriction of $v$ to $V_0$. The matrix of $v$ relative
to a direct sum decomposition $V_R=V_0\oplus V_1$ is of the form
$\pmatrix{v_0&*\cr 0&I}$, from which one gets at once $\det v_0=\det
v=\tau^{}_V(\gamma_0^{-1}\gamma,u)$.  It also follows that $\Ker v_0=\Ker v$
and that the inclusion $V_0\mono V_R$ induces an isomorphism $\Coker v_0\iso
\Coker v$,  so we deduce from (5.2)
 the exact sequence (5.3). \cqfd
\bigskip

\ind The order of $\det \gamma_0$ is $r(g-1)$ (1.9), so we can choose
$\gamma_0$ so that  $\delta=z^{-(g-1)}\gamma_0$ belongs to $SL_r(K)$.
 \medskip

 {\pc PROPOSITION} 5.4.-- {\it Let $T$ be an integral scheme, and $E$ a
vector bundle on $X\times T$, with a trivialization $\rho$ over $X^*\times
T$, such that $\wedge^r\rho$ extends to a trivialization of $\bigwedge^rE$.
Let $\mu:T\rightarrow {\cal Q}$ be the corresponding morphism {\rm (2.1)}.
Assume that for some $t\in T$, the bundle $E_{|X\times\{t\}}((g-1)p)$ has
trivial cohomology. Then the determinant bundle  $\det R\Gamma_T(E)$
{\rm (3.8)} is isomorphic to the line bundle $\mu^*{\cal L}_{\chi}$,  and
the theta divisor $\Theta_{(g-1)p}$ is the pull back of the divisor
$\delta(\div \tau^{}_V)$.}

\ind Since ${\cal L}_\chi={\cal O}_{\cal Q}\bigl(\delta(\div\tau^{}_V)\bigr)$
(4.13), the first assertion follows from the second one, which is local
over $T$. Therefore we may assume that $T=\Sp(R)$, and that
 $\mu$ is defined by an element $\gamma$ of  $SL_r\bigl(R((z))\bigr)$
(2.5). The vector bundle $E((g-1)p)$ is defined by the element $z^{g-1}\gamma$
of $GL_r\bigl(R((z))\bigr)$.  By shrinking $\Sp(R)$ if necessary, we may
assume that there exists an automorphism $u$ of $V_R$ such that $u\equiv
a(z^{g-1}\gamma_0^{-1}\gamma)$ (mod.~$\End^f(V_R)$) (4.4); the result then
follows from Cor. 5.3. \cqfd \medskip

{\pc COROLLARY} 5.5.-- {\it The pull back $\pi^*{\cal L}$ is the line bundle
${\cal L}_\chi$ on ${\cal Q}$ associated to the character $\chi$.}
\ind By prop. 6.4 below we can  write ${\cal Q}$ as a direct limit of
integral subvarieties $Q_n$.  For $n$ large enough, some points of $Q_n$ will
correspond to vector bundles $E$ on $X$   such that $E((g-1)p)$ has trivial
cohomology. Therefore by (5.4) the line bundles $\pi^*{\cal L}$ and ${\cal
L}_\chi$ have isomorphic restrictions to $Q_n$ for $n$ large enough, hence
they are isomorphic. \cqfd

\vskip 1.7cm
{\bf 6. The group ${\bf SL}_r(A_X)$}
\smallskip
\ind The next sections will be devoted to descend from the ind-variety
${\cal Q}$ to its quotient ${\bf SL}_r(A_X)\bk{\cal Q}$, which is isomorphic
to the moduli stack ${\cal SL}_X(r)$. In order to do this we will need an
important technical property of the ind-varieties ${\cal Q}$ and ${\bf
SL}_r(A_X)$, namely that they are {\it integral}. We first study a particular
case (from which we will deduce the general case): the group ${\bf
SL}_r(A_X)$ when $X={\bf P}^1$ and $p=0$. This is simply
 the $k$-group
${\bf SL}_r(k[t])$ with  $t=z^{-1}$.
\medskip
{\pc PROPOSITION} 6.1.-- {\it  The $k$-group ${\bf SL}_r(k[t])$ is the direct
limit of an
increasing sequence $(\Gamma^{(N)})_{N\ge 1}$ of subvarieties  which are
integral, normal, and locally complete intersections.}

 \ind    For any \al a $R$, define  $\Gamma^{(N)}(R)$ as the set of
matrices of degree $\le N$ in $SL_r(R[t])$. The $k$-space $\Gamma^{(N)}$ is
represented by a closed subvariety of $M_r(k)^{N+1}$,
 defined by the equation $\displaystyle \det(\sum_{n=0}^N A_n t^n)=1$. In
other words, $\Gamma^{(N)}$ is the fibre over $1$ of the
map $\det : M_r(k)^{N+1}\longrightarrow S_{rN}$ (we denote by $S_d$  the
space of polynomials in $t$ of degree $\le d$).

\ind Let $\Gamma^{(N)}_0$ be the open subset of $\Gamma^{(N)}$ consisting
of matrices $\displaystyle A(t)=\sum_{n=0}^N A_n t^n$ with $\rk(A_N)=N-1$
(the equality $\det A(t)=1$ forces $\rk(A_N)\le N-1$). Let us first prove
that {\it the map $\det : M_r(k)^{N+1}\longrightarrow S_{rN}$ is
 smooth along} $\Gamma^{(N)}_0$. Let $A(t)\in \Gamma^{(N)}_0$, and let
$M(t)=A(t)^{-1}$. The differential of $\det$ at $A(t)$ is the map
$\displaystyle B(t)\mapsto\Tr M(t)B(t)$.
The minor $M_{ij}(t)$ is of degree $\le (r-1)N$, and its highest degree
coefficient is the corresponding minor for the matrix $A_N$.
Since $A(t)$ belongs to $\Gamma^{(N)}_0$ there exist indices $i,j$ such
that $M_{ij}(t)$ has degree exactly $(r-1)N$.  Then the  minors
$M_{i1}(t),\ldots,M_{ir}(t)$, viewed as elements of  $H^0({\bf P}^1,{\cal
O}_{{\bf P}^1}((r-1)N))$, have no common zeros:  this is clear  at infinity,
because $M_{ij}(t)$ does not vanish, and on the affine line it follows from
the formula $\sum_kM_{ik}(t)A_{ki}(t)=1$. Therefore the usual resolution
for the ideal spanned by the maximal minors of a  matrix of type $(r,r-1)$
gives an exact sequence $$0\rightarrow {\cal O}_{{\bf P}^1}^{r-1}\
\hfl{A[i]}{}\  {\cal O}_{{\bf P}^1}(N)^r \ \ghfl{(M_{i1},\ldots,M_{ir})}{}\
{\cal O}_{{\bf P}^1}(rN)\rightarrow 0\ ,$$
 where $A[i]$ is obtained by deleting the $i$-th column of
$A(t)$  (see e.g. [P-S], lemme 3.1). Taking cohomology we see that the map
$(P_j)\mapsto\sum_kP_kM_{ik}$ of  $S_N^r$ into $S_{rN}$ is
surjective, which implies the surjectivity of $T_{A(t)}(\det)$. Therefore
$\Gamma^{(N)}_0$ is smooth, with the expected dimension $r^2(N+1)-(rN+1)$.

\ind We will now prove that $\Gamma^{(N)}_0$ is {\it irreducible}, and
that  $\Gamma^{(N)}\moins \Gamma^{(N)}_0$ {\it is of codimension $\ge 2$}. An
element $A(t)$ of $\Gamma^{(N)}$ can be viewed as an homomorphism $A(t):{\cal
O}_{{\bf P}^1}^r\longrightarrow {\cal O}_{{\bf P}^1}(N)^r$, which is
bijective over ${\bf P}^1\moins\{\infty\}$. Let us denote as before by
${\cal O}=k[[z]]$  the complete local ring of ${\bf P}^1$ at $\infty$. Then
the   cokernel of $A(t)$ is of the form ${\cal O}/(z^{d_1})\oplus
\ldots\oplus {\cal O}/(z^{d_r})$, with  $0\le d_1\le \ldots\le d_r$ and $\sum
d_i=rN$. The elements $z^{d_1},\ldots,z^{d_r}$ are the invariant factors of
the matrix $A(t)$ at $\infty$ ({\it i.e.} of the matrix $z^N A(z^{-1})$ over
the ring $k[[z]]$). In particular the case $(0,\ldots,0,rN)$ corresponds
exactly to $\Gamma^{(N)}_0$.

\ind Let ${\bf d}=(d_1,\ldots,d_r)$ be a sequence of integers satisfying
the above properties. Let us denote by $C_{\bf d}$ the ${\cal O}$-module
$\oplus{\cal O}/(z^{d_i})$. Using the local coordinate $z$ we can identify
the  $k$-vector space $\Hom({\cal O}_{{\bf P}^1}(N)^r,C_{\bf d})$ with
$C_{\bf d}^r$.
 Let $H_{\bf d}$ be the open subset
of this vector space consisting of those homomorphisms $\varphi$ such that
$\varphi(-1): H^0({\bf P}^1,{\cal O}_{{\bf P}^1}(N-1)^r)\longrightarrow
C_{\bf d}$ is bijective. This means that the vector bundle $\Ker \varphi$ is
trivial; it admits a unique trivialization $\tau$ such that the composite map
$A_\varphi(t):{\cal O}_{{\bf P}^1}^r\phfl{\tau}{} \Ker \varphi\mono{\cal
O}_{{\bf P}^1}(N)^r$ is the identity at $0$. Let us  consider the map
$p_{\bf d}:H_{\bf d}\times SL_r(k)\longrightarrow \Gamma^{(N)}$ defined by
$p_{\bf d}(\varphi,B)=A_\varphi(t)B$. The image of $p_{\bf d}$ is the locally
closed subvariety of  $\Gamma^{(N)}$ consisting of matrices $A(t)$ with
invariant factors at infinity $(z^{d_1},\ldots,z^{d_r})$. We see that these
subvarieties are irreducible; in particular, {\it the open subset
  $\Gamma^{(N)}_0$ is irreducible.}

\ind The automorphism group $G_{\bf d}:=\Aut_{\cal O}(C_{\bf d})$ acts freely
on $H_{\bf d}$, and $p$ clearly  factors through this action. The group
$G_{\bf d}$ is an affine algebraic group, which can be realized as an open
subvariety of the space $\displaystyle \bigoplus_{i,j} \Hom^{}_{\cal O}
\bigl( {\cal O}/(z^{d_i}),{\cal O}/(z^{d_j})\bigr)$. An easy computation
gives $\dim G_{\bf d}=(2r-1)d_1+(2r-3)d_2+\ldots +d_r\ge rN+2$ if
$d_{r-1}\not=0$. Since $\dim H_{\bf d}= r^2N$, we conclude that    $$\dim
\Gamma^{(N)}\moins \Gamma^{(N)}_0 \le (r^2N+r^2-1)-(rN+2)=\dim
\Gamma^{(N)}_0-2\ .$$ \ind Since  $\Gamma^{(N)}$ is  defined by $rN+1$
equations in $k^{r^2(N+1)}$, every component of  $\Gamma^{(N)}$ has dimension
$\ge r^2(N+1)-(rN+1)=\dim \Gamma^{(N)}_0$. We conclude that $\Gamma^{(N)}$ is
irreducible, and is a (global) complete intersection in $M_r(k)^{N+1}$. In
particular it is locally complete intersection, hence Cohen-Macaulay, and
normal by Serre's criterion. \cqfd  \bigskip {\it Remark} 6.2.-- Let
$\Gamma^{(N)}_{\bf d}$ be the image of $p_{\bf d}$; it follows easily from
the proof that $p_{\bf d}$ induces an isomorphism of $(H_{\bf d}/G_{\bf
d})\times SL_r(k)$ onto   $\Gamma^{(N)}_{\bf d}$. So we get a stratification
of $\Gamma^{(N)}$ by the smooth subvarieties $\Gamma^{(N)}_{\bf d}$, which
admit a very explicit description. One sees easily for instance that the
variety $\Gamma^{(N)}$ is rational. \bigskip \ind (6.3) We now come back to
the general case. Let us say that an ind-scheme is {\it reduced} (resp. {\it
irreducible}, resp. {\it integral}) if it is the direct limit of an
increasing sequence of reduced (resp. irreducible, resp. integral) schemes.
\medskip

{\it Lemma} 6.3.-- {\it Let $X$ be an ind-scheme, limit of an increasing
sequence of  schemes.}

\ind a) {\it If $X$ is reduced, and is the direct limit of an increasing
sequence $(X_n)$ of
 schemes, then} $X=\limind (X_n)_{\rm red}$.

\ind b) {\it If $X$ is covered by reduced open sub-ind-schemes, $X$ is
reduced.}

\ind c) {\it $X$ is integral if and only if it is reduced and irreducible.}
\ind d) {\it Let $V$ be a scheme. If $V\times X$ is
 integral), $X$ is  integral. }

\ind Let us prove a). Let $(Y_n)$ be an increasing sequence of reduced
schemes such that $X=\limind Y_n$. We have to prove that any morphism $f$
from an affine scheme into $X$ factors through some $(X_n)_{\rm red}$. But
$f$ factors through some $Y_p$, and the inclusion $Y_p\mono X$ factors
through some $X_q$. Since $Y_p$ is reduced, $f$ factors through $(X_q)_{\rm
red}$.

\ind Let us prove b). Write $X=\limind X_n$; we want to show that given
$p\in {\bf N}$, the inclusion $X_p\mono X_n$ factors through $(X_n)_{\rm
red}$ for $n$ large enough. Since $X_p$ is quasi-compact it is enough to
prove this statement locally over $X_p$, so we are reduced to the case where
$X$ is reduced; then it follows from a).

\ind The assertion c) follows from a). Let us prove d). Let $(T_n)$ be an
increasing sequence of reduced schemes such that $V\times X = \limind T_n$.
Let $p:V\times X\rightarrow X$ denote the second projection. Choose a point
$v\in V(k)$, and let $s_v:X\rightarrow V\times X$ be the section of $p$
defined by $s_v(y)=(v,y)$.  Since $X$ is an ind-variety, the induced morphism
$p:T_n\rightarrow X$ factors through a subvariety $T'_n$ of $X$, which we may
assume to be reduced (resp. irreducible) if $T_n$ is reduced (resp.
irreducible). Let $S$ be an affine scheme, and $f:S\rightarrow X$ a morphism;
writing $f=p\rond s_v\rond f$ we see that $f$ factors through $T'_n$ for some
$n$. Therefore  $X$ is the direct limit of the  varieties $T'_n$. \cqfd
\vskip 1cm

{\pc PROPOSITION} 6.4.-- {\it The ind-varieties ${\cal Q}$ and
${\bf SL}_r(A_X)$ are integral.}

\ind The ind-variety  ${\cal Q}$ is reduced by  thm.
2.5, prop. 6.1 and lemma 6.3 b), and irreducible by  prop. 2.6 e). To prove
the result for ${\bf SL}_r(A_X)$, we'll use the well-known fact that  the
open substack ${\cal SL}_X(r)^{s}$ of ${\cal SL}_X(r)$ parametrizing stable
bundles is the quotient  of a smooth variety $H$ by a group $GL_M(k)$ (see \S
8 for an explicit construction). Consider the cartesian diagram
$$\diagram{
H'&\hfl{f'}{} &{\cal Q}\cr
\vfl{\pi'}{}&&\vfl{}{\pi}\cr
H &\hfl{f}{} &{\cal SL}_X(r) \cr}$$\vskip-10pt
with $H'=H\times_{{\cal SL}_X(r)}{\cal Q}$. Reducing $H$ if necessary  we
may assume that $H'$ is isomorphic to ${\bf SL}_r(A_X)\times H$ (prop. 3.4).
The ind-variety ${\cal Q}$ is integral and the morphism $f'$ is smooth with
connected fibres (it makes $H'$ a $GL_M(k)$-torsor over ${\cal Q}$);
therefore  $H'$ is integral, and so is ${\bf SL}_r(A_X)$ by lemma 6.3 d).
\cqfd \bigskip
{\pc COROLLARY} 6.5.-- {\it ${\cal Q}$ is the direct limit of the integral
projective varieties ${\cal Q}^{(N)}_{\rm red}$.}
\ind This follows from  prop. 6.4 and lemma 6.3, a). \cqfd
\bigskip

{\pc COROLLARY} 6.6.-- {\it Every character $\chi:{\bf
SL}_r(A_X)\longrightarrow {\bf G}_m$ is trivial.} \ind We claim that the
derivative of $\chi$ (considered as a function  on ${\bf SL}_r(A_X)$) is
everywhere $0$. In fact, since $\chi$ is a homomorphism, this is equivalent
to saying that the Lie algebra homomorphism $L(\chi):{\goth
sl}_r(A_X)\longrightarrow k$ is zero. But for any commutative ring $R$ the
Lie algebra ${\goth sl}_r(R)$
 is equal to its commutator
algebra,  so any Lie algebra homomorphism of ${\goth
sl}_r(A_X)$ into $k$ is trivial.
\ind Let us write ${\bf SL}_r(A_X)$ as the limit of a sequence $V_n$ of
integral varieties. The restriction of $\chi$ to $V_n$ has again zero
derivative, hence is constant. Since $1$ belongs to $V_n$ for $n$
large enough, one has $\chi^{}_{|V_n}=1$ for all $n$, that is $\chi=1$. \cqfd

\bigskip
\ind This has the following interesting consequence:
\smallskip
{\pc PROPOSITION} 6.7.-- {\it There is a unique embedding  of  ${\bf
SL}_r(A_X)$ in $\widehat{\bf SL}_r(K)$ lifting the inclusion ${\bf
SL}_r(A_X)\i {\bf SL}_r(K)$. The corresponding
 embedding $i:{\goth sl}_r(A_X)\mono  \widehat{\goth sl}_r(K)$ is given in
terms of the decomposition $\widehat{\goth sl}_r(K)={\goth sl}_r(K)\oplus k$
{\rm (4.10)} by $i(\alpha)=(\alpha,0)$.}

\smallskip
\ind The unicity of the lifting follows from 6.6. To prove the existence,
choose an element $\delta$ in $SL_r(K)$ such that the
  bundle $E_{\delta}((g-1)p)$
has trivial cohomology. The subspace $V:=z^{-(g-1)}\delta^{-1}(A_X^r)$ is a
supplement of ${\cal O}^r$ in $K^r$ (5.2), and the elements of ${\bf
SL}_r(A_X)$ preserve $\delta(V)$, so the universal extension splits over
${\bf SL}_r(A_X)$ (4.7). To prove the assertion about the Lie algebras, the
simplest way is to notice that  the embedding of ${\goth sl}_r(A_X)$ in
$\widehat{\goth sl}_r(K)$ will  also be unique, so it is enough to check that
$i$ is a Lie algebra homomorphism. By prop. 4.10, we must prove
$\Res_0\Tr({d\alpha\over dz}\,\beta)=0$ for $\alpha,\beta$ in ${\goth
sl}_r(A_X)$; but this is a consequence of the residue theorem. \cqfd

\vskip 1.7cm
{\bf 7. The space of sections of the determinant bundle}
\smallskip

\ind (7.1) The aim of this section is to identify the space of sections of
the (powers of) the determinant bundle over the moduli stack ${\cal SL}_X(r)$
in group-theoretic terms. We first need some general formalism about descent.
We will consider   a $k$-space $Q$ and a $k$-group $\Gamma$ acting on $Q$ (on
the left). This means that we are given a morphism $m:\Gamma\times
Q\longrightarrow Q$ satisfying the usual conditions of a group action. Let
$\Gamma\bk Q$ be the quotient stack (3.2), and $\pi:Q\rightarrow \Gamma\bk Q$
be the canonical map.   \ind We suppose given a line bundle ${\cal M}_0$ on
$\Gamma\bk Q$ (3.7). Its pull back ${\cal M}=\pi^*{\cal M}_0$ to $Q$
 has a canonical {\it
$\Gamma$-linearization}, that is an isomorphism  $\varphi:m^*{\cal M}\iso
pr_1^*{\cal M}$ satisfying the usual cocycle condition. Though we will not
need it, let us observe that conversely, any line bundle on $Q$ with a
$\Gamma$-linearization comes by descent theory from a uniquely determined
line bundle on $\Gamma\bk Q$.

\ind We'll say that a section $s\in
H^0(Q,{\cal M})$  is {\it $\Gamma$-invariant} if $\varphi(m^*s)= pr_1^*s$.
We will need the following formal lemma about quotient stacks:
\medskip
{\it Lemma} 7.2.-- {\it The map $\pi^*: H^0(\Gamma\bk Q, {\cal M}_0)
\longrightarrow  H^0(Q, {\cal M})$ induces an isomorphism of $H^0(\Gamma\bk
Q, {\cal M}_0)$ onto the space of $\Gamma$-invariant sections of ${\cal M}$.}

 \ind  Since $\pi\rond m=\pi\rond q$, the
pull back of a section of ${\cal M}_0$  is $\Gamma$-invariant.
Conversely, let $s$ be a $\Gamma$-invariant section in $H^0(Q, {\cal M})$,
and let  $\mu$ be a morphism of a scheme $S$ into $\Gamma\bk Q$. Recall (3.2)
that $\mu$ corresponds to a diagram

\vskip -10pt
$$\diagram{
P&\hfl{\tilde\mu}{}&Q\cr
\vfl{}{}&&\vfl{}{}\cr
S&\hfl{\mu}{}&\Gamma\bk Q}$$\vskip -10pt
where $P$ is a $\Gamma$-torsor over $S$ and the map $\tilde\mu$ is
$\Gamma$-equivariant. By construction the section $\tilde \mu^*s$ over $P$ is
$\Gamma$-invariant (in the preceding sense); we want to show that it is the
pull back of a unique section $s_\mu$ over $S$. By standard descent theory,
it is enough to check this locally for the faithfully flat topology, so we
can suppose
 $P=\Gamma\times S$. Saying that $\tilde \mu^*s$  is
$\Gamma$-invariant means that for any map $\nu:T\rightarrow \Gamma$, where
$T$ is a scheme, the section $(\nu\times 1_S)^*\tilde \alpha^*s$ on  $T\times
S$  satisfies the usual descent condition with respect to the projection
$T\times S\rightarrow S$. Therefore this section descends to a unique section
$s_\mu\in H^0(S,\mu^*{\cal M}_0)$, which  is clearly independent of $T$, and
satisfies the required property. \cqfd \bigskip

\ind (7.3) In this situation, each element of $\Gamma(R)$ gives an
automorphism of the $k$-space $Q_R:=Q\times \Sp(R)$, hence acts on the space
$H^0(Q_R,{\cal M}_R)$; we get in this way a representation of the (abstract)
group $\Gamma(k)$ in the space $V:=H^0(Q,{\cal M})$. If $Q$ is a scheme, the
space $H^0(Q_R,{\cal M}_R)$ is canonically isomorphic to $V_R$, so the above
representation is algebraic in the sense that it is given by a morphism of
$k$-groups $\Gamma\longrightarrow \Aut(V)$. This is no longer true when $Q$
is only an ind-scheme, because inverse limits do not commute with tensor
products. They do however commute with tensor products by finite-dimensional
algebras over $k$, so what we get is a morphism
 of $\Gamma$ into $\Aut(V)$ viewed as  functors on the category of
finite-dimensional \al as. In particular the homomorphism
$\Gamma(k[\varepsilon])\longrightarrow
\Aut_{k[\varepsilon]}(V\otimes_kk[\varepsilon])$ defines in the usual way a
representation of the Lie algebra  $\Lie(\Gamma)$ on $V$. \medskip
{\pc PROPOSITION} 7.4.-- {\it Suppose $\Gamma$ and $Q$ are integral
ind-varieties {\rm (6.3)}. Let \break $s\in V=H^0(Q,{\cal M})$. The following
properties are equivalent:}

\ind (i) {\it The section $s$ is $\Gamma$-invariant;}

\ind (ii) {\it The element $s$ of $V$ is invariant under the action of
$\Gamma(k)$};

\ind (iii) {\it $s$ is annihilated by $\Lie(\Gamma)$.}

\ind (i) implies that for every \al a $R$ the image of $s$ in
$H^0(Q_R,{\cal M}_R)$ is invariant under $\Gamma(R)$; taking $R=k$ (resp.
$R=k[\varepsilon]$) gives (i) (resp. (ii)). \ind Suppose (ii) holds. Then the
section $\sigma=\varphi(m^*s)-pr_1^*s$ on $\Gamma\times Q$ vanishes by
restriction  to $\{\gamma\}\times Q$ for all $\gamma\in \Gamma(k)$; in
particular, its value at any $k$-point of $\Gamma\times Q$ is zero. Since
$\Gamma\times Q$ is reduced, this implies $\sigma=0$, hence (i).

\ind Suppose (iii) holds. Let $q\in Q(k)$, and let $i_q:\Gamma\mono
\Gamma\times Q$ denote the injection $\gamma\mapsto(\gamma,q)$. The line
bundle $i_q^*{\cal M}$  is trivialized once we choose a generator of ${\cal
M}$ at $q$,  so  we may consider $i_q^*\sigma$ as a function over $\Gamma$:
its value at a point $\gamma\in\Gamma(R)$ is obtained by evaluating  the
section $\gamma^*s-s$ at $q$. The hypothesis (iii) means that the derivative
of this function is identically zero. As in the proof of cor. 6.6 this
implies that $\sigma$ vanishes on $\Gamma\times\{q\}$ for all $q$ in $Q(k)$,
which implies as before $\sigma=0$. \cqfd  \bigskip

\ind (7.5) Let  $G$ be a $k$-group, $H$ a subgroup of $G$, and $Q$ the
quotient $G/H$. The group $G$ acts on $Q$ by left multiplication. Recall that
we have  associated to each character $\chi$ of $H$ a line bundle ${\cal
L}_\chi$ on $Q$ (3.9). We claim that this line bundle admits a {\it canonical
$G$-linearization}. The easiest way to see that is to notice that the
quotient stack $G\bk Q$ can be canonically identified with the classifying
space $BH$, with the morphism $Q\rightarrow BH$ induced by the structural map
$G\rightarrow \Sp(k)$. We have seen in (3.9) that ${\cal L}_\chi$ is the pull
back of a line bundle on $BH$, hence our assertion.  \smallskip

\ind (7.6) Let us now go back to our situation, and consider the action of
the ind-group  $\widehat{\bf SL}_r(K)$ on the ind-variety ${\cal Q}$.
According to (7.5) the line bundle ${\cal L}_\chi$ on ${\cal Q}$ admits a
canonical
 $\widehat{\bf SL}_r(K)$-linearization.
  We therefore obtain an action of the Lie algebra
$\widehat{\goth sl}_r(K)$ on the space $H^0({\cal Q}, {\cal L}_\chi)$, and
similarly on the spaces  $H^0({\cal Q}, {\cal L}_{\chi}^c)$ for all $c\in{\bf
N}$. The identification of this representation is an important result of
Kumar and Mathieu [Ku,M]. Before stating it, we need to recall a few facts
about representations of  Kac-Moody algebras, for which we refer to [K].
\ind Let us introduce the triangular decomposition ${\goth sl}_r(K)={\goth
n}_-\oplus{\goth h}\oplus{\goth n}_+$, where
 ${\goth h}$ is the Cartan algebra of diagonal matrices in ${\goth sl}_r(k)$
and ${\goth n}_+$ (resp. ${\goth n}_-$) the Lie subalgebra of matrices
$\sum_{n\ge 0}A_nz^n$ (resp. $\sum_{n\ge 0}A_nz^{-n}$) such that the matrix
$A_0$ is strictly  upper (resp. lower) triangular. The most interesting
representations of the Lie algebra $\widehat{\goth sl}_r(K)$ are the
so-called {\it integrable highest weight representations}\footnote{\parindent
0.5cm $^1$}{\eightpoint In [K] they are defined as  representations of
$\widehat{\goth sl}_r(k[z,z^{-1}])$, but we will see in (7.10) below  that
they extend to Laurent series.}\parindent0cm; they are associated to a
dominant weight $\lambda$ of the simple Lie algebra ${\goth sl}_r(k)$ and  an
integer $c\ge \langle\lambda,\tilde \alpha\rangle$, where $\tilde \alpha$ is
the highest root  of ${\goth sl}_r(k)$.  The highest weight  representation
$V_{\lambda,c}$ corresponding to the weight $\lambda$ and the integer $c$ is
characterized by the  following properties ([K], 9.10): it is irreducible, an
element $t$ of the central factor $k\i \widehat{\goth sl}_r(K)$ acts as the
homothety of ratio $ct$, and there exists a vector $v\in V_{\lambda,c}$ which
is annihilated by ${\goth n}_+$ and satisfies $H\,v=\lambda(H)v$ for all $H$
in ${\goth h}$. The  vector $v$, which is  uniquely determined up to a
scalar, is called a highest weight vector of the representation.

  \ind We will be mainly interested in the case
$\lambda=0$;  the corresponding representation $V_c$ ($c\ge 0$) is sometimes
called the {\it basic representation of level} (or charge) $c$. In this case
the annihilator of a highest weight vector $v_c\in V_c$ is ${\goth
sl}_r({\cal O})$.

 \medskip
{\pc THEOREM} 7.7 (S. Kumar, Mathieu).-- {\it The space $H^0({\cal Q}, {\cal
L}_\chi^c)$ is isomorphic (as a  $\widehat{\goth sl}_r(K)$-module) to the
dual of the basic representation $V_c$ of level $c$ of  $\widehat{\goth
sl}_r(K)$.} \ind This theorem is proved in [Ku] and [M], with one important
difference. Both S. Kumar and Mathieu define the structure of ind-variety on
$SL_r(K)/SL_r({\cal O})$ in an {\it ad hoc} way, using representation theory
of Kac-Moody algebras; we must show that it coincides with our functorial
definition. For instance Kumar, following Slodowy [Sl], considers the
representation $V_c$ for a fixed $c$, and a highest weight vector $v_c$. The
subgroup $SL_r({\cal O})$ is the stabilizer of the line $kv_c$ in ${\bf
P}(V_c)$, so the map $g\mapsto gv_c$ induces  an injection  $i_c:
SL_r(K)/SL_r({\cal O})\mono {\bf P}(V_c)$. Let $U$ be the subgroup of
$SL_r({\cal O})$ consisting of matrices $A(z)$ such that $A(0)$ is
upper-triangular with diagonal coefficients equal to $1$; to each element $w$
of the Weyl group is associated a ``Schubert variety" $X_w$ which is a finite
union of orbits of $U$. It turns out that the image under $i_c$ of  $X_w$
 is actually contained in some
finite-dimensional projective subspace ${\bf P}_w$ of ${\bf P}(V_c)$, and is
Zariski closed in ${\bf P}_w$.   This defines on $X_w$ a structure  of
reduced projective variety, and a structure of ind-variety on
$SL_r(K)/SL_r({\cal O})=\limind X_w$.

\ind To check that this ind-variety coincides with ${\cal Q}$, we  will use
the fact that the map $i_c$ is actually a morphism of ind-schemes of ${\cal
Q}$ into ${\bf P}(V_c)$ (which is the direct limit of its finite-dimensional
subspaces). In fact, we will prove in the Appendix below, following
G.~Faltings, that  the integrable representation $V_c$ of $\widehat{\goth
sl}_r(K)$ can be ``integrated" to an {\it algebraic} projective
representation of $\widehat{\bf SL}_r(K)$, that is a morphism of $k$-groups
$\widehat{\bf SL}_r(K)\longrightarrow PGL(V_c)$.  We claim that $i_c$ is an
{\it embedding}. It is injective by what we said above; let us check that it
induces an injective map on the tangent spaces. Since it is equivariant under
the action of $\widehat{\bf SL}_r(K)$ it is enough to prove this at  the
origin $\omega$ of ${\cal Q}$; then it follows from the fact  that the
annihilator of $v_c$ in the Lie algebra  $\widehat{\goth sl}_r(K)$ is ${\goth
sl}_r({\cal O})$.

\ind Therefore the restriction of $i_c$ to each of the subvarieties  ${\cal
Q}^{(N)}$ is  proper, injective, and injective on the tangent spaces, hence
is an embedding (in some finite-dimensional projective subspace of ${\bf
P}(V_c)$). Each $X_w$ is contained in some ${\cal Q}^{(N)}$, and therefore is
a closed subvariety of ${\cal Q}^{(N)}_{\rm red}$. Each orbit of $U$ is
contained in some $X_w$; since the $X_w$'s define an ind-structure, each
${\cal Q}^{(N)}$ is contained in some $X_w$, so that ${\cal Q}^{(N)}_{\rm
red}$ is a subvariety of $X_w$. Since ${\cal Q}$ is the direct limit of the
${\cal Q}^{(N)}_{\rm red}$, the two ind-structures coincide, and the theorem
follows from the Kumar-Mathieu theorem.~\cqfd \bigskip

\ind It remains to descend to the quotient $\Gamma\bk{\cal Q}$, i.e. to
apply  prop. 7.4 in the case where $Q$ is the ind-variety ${\cal Q}$, and
$\Gamma$ the ind-group ${\bf SL}_r(A_X)$. The quotient stack $\Gamma\bk{\cal
Q}$ is the moduli stack  ${\cal SL}_X(r)$ (3.4); we take for ${\cal M}_0$ a
power ${\cal L}^c$ $(c\in{\bf N})$ of the determinant bundle ${\cal L}$  on
${\cal SL}_X(r)$. What corresponds to ${\cal M}$ is the line bundle ${\cal
L}_\chi^c$ on ${\cal Q}$ (cor. 5.5).  Since ${\cal L}_\chi$ is the pull back
of ${\cal L}$, it has a canonical ${\bf SL}_r(A_X)$-linearization; on the
other hand, it has a natural $\widehat{\bf SL}_r(K)$-linearization (7.5), and
we know that the inclusion of ${\bf SL}_r(A_X)$ in ${\bf SL}_r(K)$ lifts
canonically to an embedding of ${\bf SL}_r(A_X)$ in $\widehat{\bf SL}_r(K)$
(6.7), which gives another ${\bf SL}_r(A_X)$-linearization of ${\cal
L}_\chi$. We claim that {\it these two linearizations are the same.} Actually
there is no choice: \medskip

{\it Lemma} 7.8.-- {\it The line bundle ${\cal L}_\chi$ admits a unique ${\bf
SL}_r(A_X)$-linearization.}

\ind Let us write $\Gamma={\bf SL}_r(A_X)$. Two $\Gamma$-linearizations
differ by an automorphism of $p^*{\cal L}_\chi$, i.e. by an invertible
function  on $\Gamma\times {\cal Q}$. Since ${\cal Q}$ is the direct limit of
the  integral projective varieties ${\cal Q}^{(N)}_{\rm red}$  (6.5), this
function is the pull back of an invertible function $f$ on $\Gamma$; the
cocycle conditions on the linearizations imply that  $f$ is a character,
hence $f=1$ by (6.6). \cqfd \bigskip \ind Therefore the action of the Lie
algebra ${\goth sl}_r(A_X)$ on $H^0({\cal Q},{\cal L}_\chi^c)$ is the
restriction via the natural embedding (6.7) of the action of $\widehat{\goth
sl}_r(K)$ on $H^0({\cal Q},{\cal L}_\chi^c)\cong V_c^*$ (thm. 7.7). Since the
ind-varieties ${\cal Q}$ and ${\bf SL}_r(A_X)$ are integral (prop. 6.4), we
can apply  lemma 7.2 and prop. 7.4, and we get:
\medskip

{\pc THEOREM} 7.9.-- {\it The space $H^0({\cal SL}_X(r),{\cal L}^c)$ is
canonically isomorphic to the space of conformal blocks $B_c(r)$, that is
the subspace  of $V_c^*$ annihilated by the Lie algebra ${\goth sl}_r(A_X)$.}
\cqfd \bigskip
{\it Example} 7.10.-- The only case where a direct computation seems
possible is the case $g=0$. We take as before $X={\bf P}^1$, $p=0$, so that
$A_{{\bf P}^1}=k[z^{-1}]$. The space $B_c(r)$ is the dual of $V_c/{\cal
U}^+V_c$, where ${\cal U}^+$ is the augmentation ideal of the envelopping
algebra ${\cal U}$ of ${\goth sl}_r(k[z^{-1}])$. By definition of a highest
weight module, $V_c$ is generated as a ${\cal U}$-module by a highest weight
vector $v_c$, and one has $V_c=kv_c\oplus{\cal U}^+v_c=kv_c\oplus{\cal
U}^+V_c$. We conclude that {\it the space $H^0({\cal SL}_{{\bf P}^1}(r),{\cal
L}^c)$ is one-dimensional for all $c$.}
\bigskip
{\it Remark} 7.11.-- One can deduce from the results of Mathieu  that the
Picard group of ${\cal Q}$ is generated by the line bundle ${\cal L}_\chi$
(see [M], prop. 15). It then follows from lemma 7.8 that {\it the Picard
group of ${\cal SL}_X(r)$ is generated by the determinant bundle ${\cal L}$}
(the corresponding statement for the moduli space ${\cal SU}_X(r)$ is proved
in [D-N]).

\vskip 1.7cm {\it Appendix to {\rm \S 7:} Integration of
integrable highest weight modules} (according to Faltings)
\ind In this appendix we want to show that integrable highest weight
representations of the Lie algebra $\widehat{\goth sl}_r(K)$ can be
integrated to algebraic representations of the group $\widehat{\bf SL}_r(K)$.
We will actually contend ourselves with a {\it projective} representation of
this group, since this is sufficient for our purpose and that the complete
result  requires  some more work. \smallskip
\ind (A.1) Let us denote by $\widehat{\goth sl}_r(k[z,z^{-1}])$  the sub-Lie
algebra  ${\goth sl}_r(k[z,z^{-1}])\oplus k$ of
$\widehat{\goth sl}_r(K)={\goth sl}_r\bigl(k((z))\bigr)\oplus k$ (4.10). Let
$V$ be an integrable highest weight module for $\widehat{\goth
sl}_r(k[z,z^{-1}])$.  We will use the integrability property through the
following important consequence:

{\it For each vector $v$ in $V$ there exists an integer $p$ such that
$A(z)\,v=0$ for every element $\displaystyle A(z)=\sum_{n\ge p}A_nz^n$ of
${\goth sl}_r(k[z,z^{-1}])$.}

 This means that the homomorphism $\pi:\widehat{\goth sl}_r(k[z,z^{-1}])
\longrightarrow \End(V)$ is continuous when $\widehat{\goth
sl}_r(k[z,z^{-1}])$ is endowed with the $z$-adic topology, $V$ with the
discrete topology, and $\End(V)$ with the topology of pointwise convergence.
It implies that  $\pi$ extends  to a continuous homomorphism -- still denoted
by $\pi$ -- from the $z$-adic completion  $\widehat{\goth sl}_r(K)$ of
$\widehat{\goth sl}_r(k[z,z^{-1}]) $ to $\End(V)$: one has $\displaystyle
\pi\bigl(\sum_{n\ge -N}A_nz^n\bigr)=\sum_{n\ge -N}\pi(A_nz^n)$, where the
second sum is {\it locally finite}, i.e. on each
 vector, all but finitely many of the endomorphisms in the sum are zero.
More generally,  for any \al a $R$, one gets by tensor product  a homomorphism
$\pi^{}_R:\widehat{\goth sl}_r(k[z,z^{-1}])\otimes_k R\longrightarrow
\End(V_R)$, which by continuity extends to  $\widehat{\goth
sl}_r\bigl(R((z))\bigr)$ (4.12).

\ind Suppose $\pi$ is the derivative of an algebraic representation (i.e. a
morphism of $k$-groups $\widehat{\bf SL}_r(K)\longrightarrow \Aut(V)$), such
that the center of $\widehat{\bf SL}_r(K)$ acts on $V$ by homotheties. Then
we get a {\it projective representation} of ${\bf SL}_r(K)$ in $V$, that is a
homomorphism $\overline{\rho}$ of ${\bf SL}_r(K)$ into  the  quotient
$k$-group  $PGL(V):=\Aut(V)/{\bf G}_m$, whose derivative
$L(\overline{\rho}):{\goth sl}_r(K)\longrightarrow \End(V)/k1_V$ coincides
with $\pi$ up to homotheties. We claim that we can always find such a
representation: \medskip
{\pc PROPOSITION} A.2.-- {\it Let
$\pi:\widehat{\goth sl}_r(K)\longrightarrow \End(V)$ be an integrable highest
weight representation. There exists a (unique) projective representation
$\overline{\rho}: {\bf SL}_r(K)\longrightarrow PGL(V)$ whose derivative
coincides with $\pi$ up to homotheties.} \ind The proof which follows has
been shown to us by G. Faltings.  \medskip
{\it Lemma} A.3.-- {\it Let $R$ be a \al a and  $\gamma$ an  element of
$SL_r\bigl(R((z))\bigr)$. Locally over $\Sp(R)$, there exists an
automorphism $u$ of $V_R$, uniquely determined up to an invertible element of
$R$, satisfying
$$u\pi^{}_R(\alpha)u^{-1}=\pi^{}_R\bigl(\Ad(\gamma)(\alpha)\bigr)
\leqno(A.4)$$ for any}   $\alpha\in  \widehat{\goth sl}_r\bigl(R((z))\bigr)$
({\it cf.} (4.12) for the definition of the adjoint action).  \ind We'll say
for short that an automorphism $u$ satisfying the above condition is {\it
associated} to  $\gamma$.

\ind Let us show first that this lemma implies the proposition. Thanks to
the unicity property, the automorphisms $u$ associated locally to $\gamma$
glue together to define a uniquely determined element $\overline{\rho}(g)$
in $PGL(V)(R)$. Still because of the unicity property, $\overline{\rho}$ is a
homomorphism of $k$-groups of ${\bf SL}_r(K)$ into $PGL(V)$. Let $\beta\in
\widehat{\goth sl}_r(K)$; the element $\overline{\rho}(\exp
\varepsilon\beta)$ of $PGL(V_{k[\varepsilon]})$ can be written as the class
of an automorphism $I+\varepsilon u_\beta$ of $V_{k[\varepsilon]}$, where
$u_\beta$ is an endomorphism of $V$ whose class in $\End(V)/k1_V$ is
$L(\overline{\rho})(\beta)$. Formula (A.4) applied to $R=k[\varepsilon]$ and
$\gamma=\exp \varepsilon\beta$ gives
$[u_\beta,\pi(\alpha)]=[\pi(\beta),\pi(\alpha)]$
 for each $\alpha$ in  $\widehat{\goth sl}_r(K)$. Since $\pi$ is irreducible
this implies that $u_\beta$ coincides with $\pi(\beta)$ up to homotheties
([K], lemma 9.3), hence the proposition. \cqfd \medskip
\ind We will prove lemma A.3 in several steps.
\ind a) Let us prove first the unicity assertion. We just need to observe
that an endomorphism $u$ of $V_R$ which commutes with
$\pi^{}_R\bigl(\widehat{\goth sl}_r\bigl(R((z))\bigr)\bigr)$
 is a homothety: for each $k$-linear form $\varphi:R\rightarrow k$, the
endomorphism $ (1_V\otimes\varphi)\rond u$ of $V$ commutes with
$\pi\bigl(\widehat{\goth sl}_r(K)\bigr)$, hence is a homothety, from which it
follows that $u$ is a homothety.  \smallskip

\ind b) Assume that $\gamma$ is the exponential of a  matrix
$\nu\in{\goth sl}_r\bigl(R((z))\bigr)$ which is either nilpotent, or of
positive order (so that its exponential is well-defined). Then the
automorphism $\Ad(\gamma)$ of $\widehat{\goth sl}_r\bigl(R((z))\bigr)$ is the
exponential of the derivation ${\rm ad}\ \nu$. Because of the continuity
property of $\pi^{}_R$ (A.1), the series $\exp\bigl(\pi^{}_R(\nu)\bigr)$ is
locally finite, hence defines an automorphism $u$ of $V_R$;  one has
$$\pi^{}_R(\Ad(\gamma)(\alpha))=\exp\bigl({\rm ad}\
\pi^{}_R(\nu)\bigr)(\alpha)=u\pi^{}_R(\alpha)u^{-1}\ ,$$ so $u$ satisfies
(A.4). \smallskip

\ind c) Let us observe that if two elements $\gamma,\delta$ of
$SL_r\bigl(R((z))\bigr)$ have associated automorphisms $u$ and $v$, then $uv$
is associated to $\gamma\delta$. If $R$ is a field, so is $R((z))$, hence any
element of $SL_r\bigl(R((z))\bigr)$ is a product of elementary matrices
$I+\lambda E_{ij}=\exp(\lambda E_{ij})$. The result then follows from b).
\smallskip

\ind d) The exponential mapping is a bijection from the space of matrices
$\displaystyle  \sum_{n\ge 1}A_nz^n$ with zero trace onto the group of
matrices $\displaystyle B(z)=I+\sum_{n\ge 1}B_nz^n$ with  determinant $1$
(the inverse bijection is given by the logarithm), so b) gives the result for
the matrices  $B(z)$. Let now $\gamma\in SL_r(R)$; locally over $\Sp(R)$ we
can again write $\gamma$ as a product of elementary matrices, so the result
follows as in c). Finally we see that the result holds for $\gamma$ in
$SL_r(R[[z]])$.

 \ind e) Assume now that the ring $R$ is local artinian; let ${\goth m}$  be
its maximal ideal and $\kappa$ its residue field. The quotient map
$R\rightarrow \kappa$ has a section, so the group $SL_r\bigl(R((z))\bigr)$ is
a semi-direct product of $SL_r\bigl(\kappa((z))\bigr)$ by the kernel $N$ of
the map $SL_r\bigl(R((z))\bigr)\rightarrow SL_r\bigl(\kappa((z))\bigr)$. The
lemma holds for $\gamma\in SL_r\bigl(\kappa((z))\bigr)$ by c). The elements
of $N$ are  of the form $I+A(z)$, where all the coefficients of $A(z)$ belong
to ${\goth m}R((z))$; since ${\goth m}$ is nilpotent,  $I+A(z)$ is the
exponential of a nilpotent matrix, hence the lemma holds for the elements  of
$N$ by c) and therefore for all elements of $SL_r\bigl(R((z))\bigr)$.
\smallskip

\ind f) We now arrive to the heart of the proof, the case  $\gamma\in
SL_r(R[z^{-1}])$. Let us observe that in this case one can {\it normalize} the
automorphism $u$ of $V_R$ in the following way. Let ${\cal U}({\goth n}_-)$
be  the enveloping algebra of  ${\goth n}_-$ (7.6), and ${\cal U}^*({\goth
n}_-)$ its augmentation ideal. The  space $V$ is spanned as a  ${\cal
U}({\goth n}_-)$-module by a highest weight vector $v$, and the quotient $V/
{\cal U}^*({\goth n}_-)V$ is one-dimensional. Therefore the $R$-module
$V_R/{\cal U}^*({\goth n_-})V_R$ is free of rank $1$. Since $\gamma$
normalizes ${\goth n}_-\otimes R$, $u$ induces an automorphism of this
R-module, so we can choose $u$ (in a unique way) so that it induces the
identity \ mod. ${\cal U}^*({\goth n_-})V_R$.
\ind Since the group ${\bf SL}_r({\cal O}_-)$ is an ind-variety, we may assume
that the \al a $R$ is finitely generated. We will prove the existence of an
associated automorphism $u$ by induction on $\dim(R)$.

 \ind Let ${\goth p}_1,\ldots,{\goth p}_s$ be the minimal prime ideals of
$R$, and $S=R\moins \cup {\goth p}_i$. Over the artinian ring $S^{-1}R$
(isomorphic to $\prod R_{{\goth p}_i}$) we can construct by e) an
automorphism $u_S$ of $V_{S^{-1}R}$ associated to $\gamma$. Let us denote by
${\cal U}$ the envelopping algebra of $\widehat{\goth sl}_r(K)$. Observe that
(A.4) is equivalent to saying that $u$ is a {\it semi-linear} endomorphism of
the ${\cal U}$-module $V_R$ relative to the automorphism $\Ad(\gamma)$ of
${\cal U}$. Since the ${\cal U}$-module $V$ is finitely presented ([K],
10.4.6), we can find an element $f$ of $S$ such that $fu_S$ comes from a
$({\cal U}\otimes R, \Ad(\gamma))$-linear endomorphism $u_f$ of $V_R$ (chase
denominators in generators and relations).
Then the class of $u_f(v)-fv$ in $V_R/{\cal U}^*({\goth n}_-)V_R$ is
annihilated by some element $s$ of $S$, so replacing $f$ by $fs$ we may
assume $u_f(v)\equiv v$ (mod.~${\cal U}^*({\goth n_-})V_R$). Moreover we can
modify $f$ so that  $fu_S^{-1}$ also comes from an endomorphism $u'$ of
$V_R$,  such that the endomorphisms $\displaystyle {u_f\over f}$ and
$\displaystyle {u'\over f}$ of $V_{R_f}$ are inverse of each other.

\ind Since $f\in S$, one has $\dim (R/f^nR)<\dim(R)$ for each $n$, hence the
induction hypothesis provides a normalized automorphism of $V_{R/f^nR}$
satisfying (A.4). These automorphisms define an automorphism $\hat u$ of the
$f$-adic completion $\widehat V_R$ of $V_R$. On the other hand, $u_f$ extends
to an endomorphism $\hat u_f$ of $\widehat V_R$; one has $\hat u_f\equiv
f\hat u$ (mod. $f^n\widehat V_R$) for all $n$, hence $\hat u_f= f\hat u$.

\ind Unfortunately $\widehat V_R$ is bigger than $V_{\widehat
R}:=V\otimes_k \widehat R$: if we choose a basis $(e_\iota)_{\iota\in I}$ of
$V$, the elements of $\widehat V_R$ are formal sums $\sum r_\iota e_\iota$,
where for every $n\ge 0$, $f^n$ divides all but finitely many of the
$r_\iota$'s. However, since $\widehat R$ is noetherian, there exists an
integer $n$ such that  $\Ann_{\widehat R}(f^n)=\Ann_{\widehat R}(f^{n+1})$,
which implies $f^n\widehat R\,\cap\,\Ann_{\widehat R}(f)=0$; therefore {\it
an element $x$ of $\widehat V_R$ such that $fx\in V_{\widehat R}$ belongs
itself to $V_{\widehat R}$.} Coming back to our situation, we deduce from the
formula $\hat u_f= f\hat u$ that $\hat u$ induces an endomorphism $\tilde u$
of $V_{\widehat R}$; using the same construction with $u^{-1}_S$ shows that
$\hat u^{-1}$ also preserves $V_{\widehat R}$, so that $\tilde u$ is an
automorphism. It acts trivially on $V_{\widehat R}/{\cal U}^*({\goth
n_-})V_{\widehat R}$, because it does mod. $f^n$ for all $n$.

\ind  By the unicity property, the automorphisms
$\displaystyle {u_f\over f}$ of $V_{R_f}$ and $\tilde u$ of $V_{\widehat
R}$ have the same image in $\Aut(V_{{\widehat R}\otimes_R R_f})$. Since the
homomorphism $R\longrightarrow R_f\times \widehat R$ is faithfully flat, they
can be glued together to define an automorphism $u$ of $V_R$, which satisfies
(A.4) because both $\displaystyle {u_f\over f}$  and $\tilde u$ do.
\smallskip \ind g) Finally  the general case follows from lemma (4.5) and
cases c), d) and f).~\cqfd

\vskip1.7cm

{\bf 8. From the moduli stack to the moduli space}
\smallskip
\ind The last step is to compare the sections of the determinant bundle (or
of its powers) over the moduli space and over the moduli stack. {\it
Throughout this section we assume $g\ge 1$} (by example 7.10 there is
essentially nothing to say in the case $g=0$).

\ind (8.1) We first review briefly the standard construction of the moduli
space (or stack) of vector bundles. For each integer $N$, we will denote by
${\cal SL}_X(r)_{N}$
 the open substack of ${\cal
SL}_X(r)$ parametrizing vector bundles $E$ on $X$ such that  $H^1(X,E(Np))$
is $0$ and $H^0(X,E(Np))$ is generated by its global sections. Let
$h(N)=\dim H^0(X,E(Np))$ ($=r(N+1-g)$). Choosing an isomorphism
$k^{h(N)}\longrightarrow H^0(X,E(Np))$ realizes $E$ as a quotient of the
bundle  ${\cal O}_X(-Np)^{h(N)}$. The stack which  parametrizes such
quotients is represented by a  smooth scheme $K_N$. Let ${\cal E}$ be the
universal quotient bundle over $X\times K_N$, and let $q:X\times
K_N\longrightarrow  K_N$ denote the second projection. The sheaf
 $q_*\bigwedge^r{\cal E}$ is the sheaf of sections of a line bundle  on
$K_N$; let $H_N$ be the complement  of the zero section in this line bundle.
By construction $H_N$ parametrizes quotients $E$ of ${\cal O}_X(-Np)^{h(N)}$
together with a trivialization of $\bigwedge^rE$. The group $GL(h(N))$ acts
on $H_N$, and  the quotient stack is ${\cal SL}_X(r)_{N}$.
\medskip

\ind We'll denote by  ${\cal SL}_X(r)^{ss}$ the open substack of  ${\cal
SL}_X(r)$ parametrizing semi-stable bundles, and by $H_N^{ss}$ the
corresponding open subset of $H_N$. We'll assume $N\ge 2g$, which  implies
that ${\cal SL}_X(r)^{ss}$ is contained in ${\cal SL}_X(r)_{N}$.

\medskip
{\it Lemma} 8.2.-- {\it The codimension of $H_N\moins H_N^{ss}$ in $H_N$  is
at least  $2$.}

\ind For each pair of integers $(s,d)$ with $0<s<r$ and $d>0$, let us
define a stack ${\cal SL}_X^{s,d}(r)$ by associating to a \al a $R$ the
groupoid of triples $(E,F,\delta)$, where  $E$ is a rank $r$ vector bundle
over $X_R$, $F$ a rank $s$ subbundle of degree $d$, and $\delta$ a
trivialization of $\bigwedge^r E$. Forgetting $F$ gives a  morphism of stacks
of ${\cal SL}_X^{s,d}(r)$ to ${\cal SL}_X(r)$; the (reduced) substack  ${\cal
SL}_X(r)\moins {\cal SL}_X(r)^{ss}$ is the union of these images (for
variable $s,d$). According to [L], cor. 2.10, the dimension of ${\cal
SL}_X^{s,d}(r)$ is $(g-1)(r^2-1+s^2-rs)-rd$, so the codimension of its image
is at least $rd$, which is $\ge 2$. Since $H_N$ is a torsor over ${\cal
SL}_X(r)_N$ the lemma follows. \cqfd

{\pc PROPOSITION} 8.3.-- {\it For any integer $c$, the restriction map
$$H^0({\cal SL}_X(r),{\cal L}^c)\longrightarrow H^0({\cal SL}_X(r)^{ss},{\cal
L}^c)$$ is an isomorphism.}
\ind Let $N\ge 2g$. Consider the diagram
\vskip -15pt
$$\diagram{
H_N^{ss}&\mono &H_N&\cr
\vfl{}{}&&\vfl{}{}&\cr
{\cal SL}_X(r)^{ss}&\mono&{\cal SL}_X(r)_N&\ .}$$
\vskip -10pt
By lemma 7.2, the  sections of ${\cal L}^c$ over ${\cal SL}_X(r)^{ss}$
(resp.  ${\cal SL}_X(r)_N$) are the invariant sections of the pull back of
${\cal L}^c$ over $H_N^{ss}$ (resp. $H_N$). But  any section over $H_N^{ss}$
extends to $H_N$ by (8.2), so the restriction map $H^0({\cal SL}_X(r)_N,{\cal
L}^c)\longrightarrow H^0({\cal SL}_X(r)^{ss},{\cal L}^c)$ is an isomorphism
for each $N$. Since any map from a scheme to ${\cal SL}_X(r)$ factors through
${\cal SL}_X(r)_N$ for some $N$, the proposition follows. \cqfd \bigskip

\ind Let ${\cal SU}_X(r)$ be the moduli space of semi-stable rank $r$ vector
bundles on $X$ with trivial determinant (the notation is meant to remind that
these correspond to unitary representations). It is usually constructed as
the geometric invariant theory quotient of $K_N^{ss}$ (8.1) by the group
$PGL(h(N))$. We have a forgetful morphism $\varphi:{\cal
SL}_X(r)^{ss}\longrightarrow  {\cal SU}_X(r)$. It is known that the
determinant bundle ${\cal L}$ is the pull back of a line bundle (that we will
still denote by ${\cal L}$) on ${\cal SU}_X(r)$ (see
 [D-N], and [Tu] for the case $g=1$).
\medskip

 {\pc PROPOSITION} 8.4.-- {\it Let $c\in {\bf N}$. The map
$\varphi^*:H^0({\cal SU}_X(r),{\cal L}^c) \longrightarrow  H^0({\cal
SL}_X(r)^{ss},{\cal L}^c)$ is an isomorphism.}
\ind Let us choose an integer
$N\ge 2g$. We claim that both spaces can be identified with the space of
$GL(h(N))$-invariant sections of the pull back of ${\cal L}^c$ to $H_N^{ss}$.
For  $H^0({\cal SL}_X(r)^{ss},{\cal L}^c)$, this follows from lemma 7.2.

\ind Let us consider the space $H^0({\cal SU}_X(r),{\cal L}^c)$. We will
write simply $H$, $K$ and ${\cal S}$ for $H_N^{ss}$, $K_N^{ss}$ and ${\cal
SU}_X(r)$. By definition of the GIT quotient, the quotient map
$p:K\longrightarrow {\cal S}$ is affine and the sheaf  ${\cal O}_{\cal S}$
is  the  subsheaf of local $PGL(h(N))$-invariant sections in $p_*{\cal
O}_{K}$. On the other hand, since the map $q:H\rightarrow  K$ is a
$k^*$-fibration, the subsheaf of  local $k^*$-invariant sections of
$q_*{\cal O}_{H}$  is ${\cal O}_{K}$. Putting things together we conclude
that the sheaf of $GL(h(N))$-invariant sections of $p_*q_*{\cal O}_H$ is
${\cal O}_{\cal S}$. Therefore for any sheaf ${\cal F}$ on  ${\cal S}$ the
space $H^0({\cal S},{\cal F})$ is the space of $GL(h(N))$-invariant sections
of the pull back of ${\cal F}$ to $H$.
 \cqfd
\bigskip
\ind Putting together prop. 8.3 and 8.4 and thm. 7.9, we obtain:
\smallskip

{\pc THEOREM} 8.5.-- {\it For all $c\in {\bf N}$, the space $H^0({\cal
SU}_X(r),{\cal L}^c)$ is canonically isomorphic to the space of conformal
blocks $B_c(r)$.} \cqfd
\bigskip
\ind According to [T-U-Y], the dimension of the space $B_c(r)$ is constant
under semi-stable degenerations of the curve $X$, and can therefore be
computed for a curve $X$ which is a union of rational curves, each of them
meeting the rest of the curve in $3$ points. According to the conjecture of
Verlinde proved in [M-S], this gives the following  formula (this formulation
has been shown to us by D. Zagier): \medskip

{\pc COROLLARY} 8.6 (Verlinde formula).-- {\it One has $$\dim H^0({\cal
SU}_X(r),{\cal L}^c)=\Bigl({r\over r+c}\Bigr)^g \sum_{{S\i [1,r+c]}\atop
|S|=r}\ \prod_{ \scriptstyle s\in S\atop \scriptstyle t\notin S}\big | 2\sin
\pi{s-t\over  r+c}\,\bigr |^{g-1}\ .\quad \carre $$}

\vskip 1.7cm

{\bf 9. Arbitrary degree }
\smallskip
\ind In this last section we will extend our results to the case of vector
bundles of arbitrary degree. We fix an integer $d$, and let ${\cal
SL}_X(r,d)$ be the moduli stack parametrizing vector bundles $E$ on $X$ of
rank $r$ with an isomorphism $\delta:{\cal O}_X(dp)\iso \bigwedge^rE$. This
stack depends only on the class of $d$ mod. $r$, so we'll loose no generality
by assuming $0< d< r$. We will still use the letter ${\cal L}$ to denote the
determinant bundle on  ${\cal SL}_X(r,d)$ (3.8).

\ind  Recall that
the {\it fundamental weights} $\varpi_1,\ldots,\varpi_{r-1}$ of ${\goth
sl}_r(k)$ are the linear forms on the Cartan algebra ${\goth h}\i {\goth
sl}_r(k)$ defined by $\displaystyle \langle
\varpi_k\,,\,H\rangle=\sum_{i=1}^kH_{ii}$. Using the notation of (7.6), we
can state the main result of this section:
 \medskip
{\pc THEOREM} 9.1.-- {\it Let $0< d< r$. The space $H^0({\cal
SL}_X(r,d),{\cal L}^c)$ is canonically isomorphic to the subspace  of
$V_{c\varpi_{r-d},c}^*$ annihilated by the Lie algebra ${\goth sl}_r(A_X)$.}
\medskip \ind The proof follows the same lines as in the degree zero case.
We choose once and for all  an  element $\gamma_d$ of $GL_r(K)$ with
determinant of order $d$.
Then ${\cal SL}_X(r,d)$ can be described as the quotient
stack $\bigl(\gamma_d^{-1}{\bf SL}_r(A_X)\gamma_d\bigr)\bk {\cal Q}$ (3.6).
Let  $\pi_d:{\cal Q}\longrightarrow  {\cal SL}_X(r,d)$ be the canonical
morphism.  \medskip
{\pc PROPOSITION} 9.2.-- {\it One has $\pi_d^*{\cal L}\cong {\cal L}_\chi$.}

\ind We will reduce this assertion to the case $d=0$  by using the following
trick. The line bundle ${\cal O}_X(-dp)$ has a natural trivialization
$\rho_0$ over $X^*$; there is a unique  trivialization $\sigma_0$ of ${\cal
O}_X(-dp)$ over $D$ such that the  element $\rho_0^{-1}\sigma_0$ of $K^*$
is equal to $(\det \gamma_d)^{-1}$ (1.5).
Let $R$ be a \al a, and
 $E$  a
vector bundle on $X_R$, of rank $r$ and degree $d$, with trivializations
$\rho$ over $X^*_R$ and $\sigma$ over $D_R$, corresponding to an element
$\delta$ of $GL_{r}\bigl(R((z))\bigr)$ (prop. 1.4). Then the triple $(E\oplus
{\cal O}_X(-dp),\rho\oplus \rho_0,\sigma\oplus \sigma_0)$ corresponds to the
matrix $\pmatrix{\delta&0\cr 0&(\det \gamma_d)^{-1}}$. If $\delta=\gamma_d
\gamma$, with $\gamma\in SL_r\bigl(R((z))\bigr)$, this matrix is the product
of $\gamma'_d:=\pmatrix {\gamma_d & 0\cr 0&(\det \gamma_d)^{-1}}$ with the
matrix $t(\gamma):=\pmatrix{\gamma&0\cr 0&1}$.  We have therefore obtained a
commutative diagram (we use a $'$ when we replace $r$ by  $r+1$ in the
objects defined in \S 3 and 4): \vskip -15pt $$\diagram{
{\cal Q}&\hfl{f}{}&{\cal Q}'\cr
\vfl{\pi_d}{}&&\vfl{}{\pi'}\cr
{\cal SL}_X(r,d)&\hfl{s}{}&{\cal SL}_X(r+1)
}$$\vskip -10pt
where $f$ is induced by the map $\gamma\mapsto \gamma'_d\,t(\gamma)$ from
${\bf SL}_r(K)$ into  ${\bf SL}_{r+1}(K)$, and $s$ associates to a \vb $E$ on
$X_R$ the \vb $E\oplus {\cal O}_X(-dp)$.

\ind Let us denote by $E$ and $E'$  the universal bundles on  $X\times {\cal
SL}_X(r,d)$ and $X\times {\cal SL}_X(r+1)$ respectively. By construction the
pull back of $E'$  by $1_X\times s$ is $E\oplus {\cal O}_X(-dp)$. Let $p$ be
the  projection from $X\times {\cal SL}_X(r,d)$ onto ${\cal SL}_X(r,d)$.  One
has $Rp_*\bigl(E\oplus{\cal O}_X(-dp)\bigr)\cong Rp_*(E)\oplus
Rp_*\bigl({\cal O}_X(-dp)\bigr)$ and the bundles $R^ip_*\bigl({\cal
O}_X(-dp)\bigr)$ are trivial, so we get $s^*{\cal L}'\cong \det
Rp_*\bigl(E\oplus{\cal O}_X(-dp)\bigr) \cong \det Rp_*(E)={\cal L}$.
Therefore  our assertion is equivalent to $f^*{\cal L}_{\chi'}\cong {\cal
L}_\chi$.

\ind  The group morphism $t:{\bf SL}_r(K)\longrightarrow
{\bf SL}_{r+1}(K)$ extends in a straightforward way to $\hat t:\widehat{\bf
SL}_r(K)\longrightarrow  \widehat{\bf SL}_{r+1}(K)$: from the decomposition
$K^r=V\oplus{\cal O}^r$ and an arbitrary decomposition $K=V_0\oplus{\cal O}$
one gets $K^{r+1}=(V\oplus V_0)\oplus {\cal O}^{r+1}$; then $\hat t$ is
defined by  $\hat t(\gamma,u)=\left(\pmatrix{\gamma&0\cr 0&1},
\pmatrix{u&0\cr 0&1_{V_0}}\right)$. This morphism maps $\widehat{\bf
SL}_r({\cal O})$ into $\widehat{\bf SL}_{r+1}({\cal O})$
 and satisfies $\hat
t(\gamma,a(\gamma))=\bigl(t(\gamma),a(t(\gamma))\bigr)$ for   $\gamma\in
SL_r\bigl(R((z))\bigr)$, from which one deduces $\chi'\rond \hat t=\chi$. By
(3.9) this implies that the pull back of ${\cal L}_{\chi'}$ by the morphism
${\cal Q}\rightarrow {\cal Q}'$ deduced from $\hat t$ is isomorphic to ${\cal
L}_\chi$.  Since ${\cal L}_{\chi'}$ is invariant under the action of
$SL_{r+1}(K)$ (7.5), we conclude that $f^*{\cal L}_{\chi'}$ is isomorphic to
${\cal L}_\chi$, which proves the proposition. \cqfd \bigskip

\ind It follows from the proposition that the line bundle ${\cal L}_\chi$
on ${\cal Q}$ has a  $\bigl(\gamma_d^{-1}{\bf
SL}_r(A_X)\gamma_d\bigr)$-linearization (7.1), which means that the map
$\delta\mapsto \gamma_d^{-1}\delta\gamma_d$ of ${\bf SL}_r(A_X)$ into ${\bf
SL}_r(K)$ lifts to  $\widehat{\bf SL}_r(K)$. It is not difficult to describe
explicitely this lifting (use the same trick as in the above proof), but we
need only to know the corresponding Lie algebra map  $i_d:{\goth
sl}_r(A_X)\longrightarrow \widehat{\goth sl}_r(K)$. As in (6.7) we just have
to exhibit one homomorphism which coincides with
$\alpha\mapsto\gamma_d^{-1}\alpha\gamma_d $ modulo the center of
$\widehat{\goth sl}_r(K)$.
\ind Recall (4.12) that the adjoint action of
$SL_r(K)$ on  $\widehat{\goth sl}_r(K)$ is given by
$$\Ad(\gamma)\ (\alpha,s)=\bigl(\gamma\alpha\gamma^{-1}\,,\,s+\Res_0
\Tr(\gamma^{-1}{d\gamma\over dz}\ \alpha)\bigr)\ .\leqno(9.3)$$ We observe
that the formula makes sense for $\gamma\in GL_r(K)$, and defines a group
homomorphism $GL_r(K)\longrightarrow \Aut({\goth sl}_r(K))$, that we still
denote by $\Ad$. Then the homomorphism $\Ad(\gamma_d^{-1})\rond i$, where $i$
is the canonical
 embedding of ${\goth sl}_r(A_X)$ into $\widehat{\goth sl}_r(K)$ (6.7),
satisfies the required conditions and is therefore equal to $i_d$.

\ind Let us denote by $\pi$ the homomorphism $\widehat{\goth
sl}_r(K)\longrightarrow \End(V_c)$. By prop.~7.4 and thm. 7.7, the space
$H^0({\cal SL}_X(r,d),{\cal L}^c)$ is canonically isomorphic to the space of
linear forms on $V_c$ which vanish on the image of $\pi(i_d(\alpha))$ for
every $\alpha$ in ${\goth sl}_r(A_X)$, i.e. of linear forms killed by ${\goth
sl}_r(A_X)$ acting on $V_c$ through the representation $\pi\rond
\Ad(\gamma_d^{-1})$. Therefore thm. 9.1 will be a consequence of the
following lemma:

\medskip

{\it Lemma} 9.3.-- {\it The representation $\pi\rond \Ad(\gamma_d^{-1})$  is
isomorphic to the  highest weight representation $V_{c\varpi_{r-d},c}$.}
\ind  By lemma (A.3), the representations $\pi$ and $\pi\rond \Ad(\gamma)$
are isomorphic for $\gamma\in SL_r(K)$, so the representation $\pi\rond
\Ad(\gamma_d^{-1})$ doesn't depend (up to isomorphism) of the choice of the
particular element $\gamma_d$; we choose for $\gamma_d$ the matrix \
diag$(1,\ldots,1,z,\ldots,z)$, where $z$ appears $d$ times, and denote by
$\pi_d$ the representation $\pi\rond \Ad(\gamma_d^{-1})$. Let $\alpha=\sum
A_nz^n\in{\goth sl}_r(K)$; an easy computation (using (9.3)) gives
$$\Ad(\gamma_d^{-1})(\alpha,s)=\bigl(\gamma_d^{-1}\alpha\gamma_d\,,\,
s+\langle \varpi_{r-d}\,,\,A_0\rangle\bigr)\ .$$

This implies
$\Ad(\gamma_d^{-1})({\goth n}_+)\i {\goth sl}_r({\cal O})$, and therefore
the highest weight vector $v_c$ of $V_c$ is annihilated by  $\pi_d({\goth
n}_+)$ (7.6). Let $H\in {\goth h}$, and $s\in k$; the above formula gives
$\pi_d\bigl((H,s)\bigr)v_c=c(s+\langle \varpi_{r-d}\,,\,H\rangle)v_c$.
Moreover the representation $\pi_d$ is irreducible. Therefore
 $\pi_d$ is isomorphic to the highest weight
  representation $V_{c\varpi_{r-d},c}$ (9.6). \cqfd
\bigskip
\ind Let ${\cal SU}_X(r,d)$ be the moduli space of semi-stable vector
bundles on $X$ of rank $r$ and determinant ${\cal O}_X(dp)$. As in \S 8 we
have a forgetful morphism $\varphi:{\cal SL}_X(r,d)^{ss}\longrightarrow
{\cal SU}_X(r,d)$. According to [D-N], the determinant bundle ${\cal L}$
itself does not descend in general to a line bundle on  ${\cal SU}_X(r,d)$;
the pull back of the ample generator ${\cal L}_{r,d}$ of $\Pic\bigl({\cal
SU}_X(r,d)\bigr)$ is the line bundle $\det R\Gamma^{}_{{\cal
SL}_X(r,d)}({\cal E}\otimes F)$ (3.8), where  ${\cal E}$ is the universal
bundle over $X\times {\cal SL}_X(r,d)$ and $F$ a \vb on $X$ of rank
$\displaystyle s:={r\over (r,d)}$ . By (3.8) we get $\varphi^*{\cal
L}_{r,d}\cong{\cal L}^s$. Now the proof of thm. 8.5 applies almost without
modification to this situation; the only point which requires some care is
lemma 8.2, where one gets $\codim(H_N\moins H^{ss}_N)=1$
 in the case $g=1$, $(r,d)=1$. Assuming $g\ge 2$ for simplicity, we obtain
\medskip
{\pc THEOREM} 9.4.-- {\it Assume $0< d< r$ and $g\ge 2$; let $\displaystyle
s={r\over (r,d)}$. The space \break$H^0({\cal SU}_X(r,d),{\cal L}_{r,d}^c)$
is canonically isomorphic to the subspace  of $V_{cs\varpi_{r-d},cs}^*$
annihilated by the Lie algebra ${\goth sl}_r(A_X)$.}~\cqfd \par

\baselineskip13pt

\vskip 2cm
\def\num#1{\item{\hbox to\parindent{\enskip [#1]\hfill}}}
\parindent=1.5cm

\centerline{{\bf REFERENCES}}
\bigskip
\num{A-D-K} E. {\pc ARBARELLO}, C. {\pc DE} {\pc CONCINI}, V. {\pc KAC}:
{\it The infinite wedge representation and the reprocity law for algebraic
curves.} Proc. of Symp. in Pure Math. {\bf 49}, 171-190 (1989).
\smallskip
\num B N. {\pc BOURBAKI}: {\it Alg\`ebre Commutative,} ch. 5 to 7. Masson,
Paris (1985). \smallskip
\num{D-M} P. {\pc DELIGNE}, D. {\pc MUMFORD}: {\it The irreducibility of the
space of curves of given genus.} Pub. math. IHES {\bf 36}, 75-110 (1969).
\smallskip
\num{D-N} J.M. {\pc DREZET}, M.S. {\pc NARASIMHAN}: {\it Groupe de Picard
des vari\'et\'es de modules de fibr\'es semi-stables sur les courbes
alg\'ebriques.} Invent. math. {\bf 97}, 53-94 (1989).
 \smallskip
\num F G. {\pc FALTINGS}: {\it A proof for the Verlinde formula.} Preprint
(1992). \smallskip

\num K V. {\pc KAC}: {\it Infinite dimensional Lie algebras.} Progress in
Math. {\bf 44}, Birkh\"auser, Boston (1983).
\smallskip
\num{Ku} S. {\pc KUMAR}: {\it Demazure character formula in arbitrary
Kac-Moody setting.} Invent. math. {\bf 89}, 395-423 (1987).  \smallskip

\num{K-N-R} S. {\pc KUMAR}, M.S. {\pc NARASIMHAN}, A. {\pc RAMANATHAN}:
{\it Infinite Grassmannian and moduli space of {\rm G}-bundles.} Preprint
(1993). \smallskip

 \num L G. {\pc LAUMON}: {\it Un analogue global
du c\^one nilpotent.} Duke math. J. {\bf 57}, 647-671 (1988).
\smallskip
\num{L-MB} G. {\pc LAUMON}, L. {\pc MORET-BAILLY}: {\it Champs
alg\'ebriques.} Pr\'epublication Universit\'e Paris-Sud (1992).
\smallskip
\num{M} O. {\pc MATHIEU}: {\it Formules de caract\`eres pour les alg\`ebres
de Kac-Moody g\'en\'erales.} Ast\'erisque {\bf 159-160} (1988).
\smallskip
\num{M-S} G. {\pc MOORE}, N. {\pc SEIBERG}: {\it Polynomial equations for
rational conformal field theories.} Physics Letters B {\bf 212}, 451-460
(1988).
\smallskip
\num{P-S} C. {\pc PESKINE}, L. {\pc SZPIRO}: {\it Dimension projective finie
et cohomologie locale.} Pub. math. IHES {\bf 42},
47-119 (1973).
\smallskip
\num{SGA 6} {\it Th\'eorie des intersections et th\'eor\`eme de
Riemann-Roch.} S\'eminaire de G\'eom\'etrie alg\'ebrique SGA 6, dirig\'e par
P. Berthelot, A. Grothendieck, L. Illusie. Lecture Notes {\bf 225},
Springer-Verlag (1971.) \smallskip

 \num{Sl} P. {\pc SLODOWY}: {\it On the geometry of Schubert varieties
attached to Kac-Moody Lie algebras.} Can. Math. Soc. Conf. Proc. {\bf 6},
405-442 (1984). \smallskip
\num{S-W} G. {\pc SEGAL}, G. {\pc WILSON}:  {\it Loop groups and equations
of KdV type.}  Pub. math. IHES {\bf 61},
5-65 (1985).
\smallskip

\num T J. {\pc TATE}: {\it Residues of differentials on curves.} Ann.
scient. \'Ec. Norm. Sup. {\bf 1} (4\up{\`eme} s\'erie), 149-159 (1968).
\smallskip
\num{Tu} L. {\pc TU}: {\it Semistable bundles over an elliptic
curve.} Adv. in Math., to appear. \smallskip
\num{T-U-Y} A. {\pc TSUCHIYA}, K. {\pc UENO}, Y. {\pc YAMADA}: {\it Conformal
field theory on
universal family of stable curves with gauge symmetries.} Adv. Studies in
Pure Math. {\bf 19}, 459-566 (1989).
\smallskip
 \num V E. {\pc VERLINDE}: {\it Fusion rules and modular transformations in
2d conformal field theory.} Nuclear Physics {\bf B300}  360-376 (1988).
  \smallskip
\num W E. {\pc WITTEN}: {\it Quantum field theory and the Jones polynomial.}
Commun. math. Phys. {\bf 121}, 351-399 (1989).  \vskip 1.7cm

\hfill\hbox to 5cm{\hfill A. Beauville, Y. Laszlo\hfill}
\smallskip

\hfill\hbox to 5cm{\hfill URA 752 du CNRS\hfill}
\smallskip

\hfill\hbox to 5cm{\hfill Math\'ematiques -- B\^at. 425\hfill}
\smallskip

\hfill\hbox to 5cm{\hfill Universit\'e Paris-Sud\hfill}
\smallskip

\hfill\hbox to 5cm{\hfill 91 405 {\pc ORSAY} Cedex, France\hfill}

\bye